\documentclass[11pt]{article}
\usepackage[latin1]{inputenc}
\usepackage[dvips]{graphicx}
\usepackage{latexsym}
\usepackage{amsmath,amssymb,amsfonts,amsthm}

\newcommand{\fek}{{\cal F}_{k-1}}

\newcommand{\sig}{\sigma\sigma^*}
\newcommand{\nrn}{\rightarrow+\infty}
\newcommand{\xrn}{\xrightarrow}

\newcommand{\ER}{\mathbb {R}}\newcommand{\EN}{\mathbb {N}}
\newcommand{\PE}{\mathbb {P}}
\newcommand{\ES}{\mathbb{E}}
\newcommand{\abs[1]}{\lvert#1\rvert}

\newcommand{\guil}{\textquotedblleft}
\newcommand{\psg}{\langle }
\newcommand{\psd}{\rangle }
\newtheorem{theorem}{ \textnormal{\bf{T\scriptsize{HEOREM}}}}
\newtheorem{prop}{\textnormal{\bf{P\scriptsize{ROPOSITION}}}}

\newtheorem{lemme}{\textnormal{\bf{L\scriptsize{EMMA}}}}

\theoremstyle{definition}

\theoremstyle{remark}

\newtheorem{Remarque}{\textnormal{\bf{R\scriptsize{EMARK}}}}

\author{Fabien Panloup\thanks{ \textit{Laboratoire de Probabilit\'es et Mod\`eles Al\'eatoires}
    \textit{Universit\'e Pierre et Marie Curie et C.N.R.S. UMR 7599}
    \textit{Paris, France}
       \textbf{e-mail:}panloup@ccr.jussieu.fr
            \textbf{Adresse postale:} Fabien Panloup, Universit\'e Pierre et Marie Curie, LPMA, bureau 4D1, 175, rue du Chevaleret F-75013 Paris, France
\textbf{t\'el\'ephone:} (+33) 1 44 27 85 10}}
\title{\textbf{Computation of the invariant measure for a Lévy driven SDE: Rate of convergence  }}

\begin{document}
\maketitle
\begin{abstract} We study the rate of convergence of some
recursive procedures based on some \textquotedblleft exact" or
\textquotedblleft approximate" Euler schemes which converge to the
invariant measure of an ergodic SDE driven by a Lévy process.
The main interest of this work is to compare the rates induced by
\textquotedblleft   exact" and \textquotedblleft   approximate"
Euler schemes. In our main result, we show  that replacing the small
jumps by a Brownian component in the approximate case   preserves
the rate induced by the exact Euler scheme for a large class of Lévy
processes.
\end{abstract}
\section{\large{Introduction}}
In a recent paper (see \cite{panloup}), we investigated
a family of several weighted empirical measures based on some Euler schemes with
decreasing step in order to approximate recursively the invariant distribution
$\nu$ of an ergodic jump diffusion process
$X=(X_t)_{t\ge0}$ solution to a SDE driven by a Lévy process.
 More precisely, let
$(\bar{X}_k)_{k\ge 1}$ be such an Euler scheme with
  sequence of decreasing steps $(\gamma_k)_{k\ge 1}$ and let $(\eta_k)_{k\ge 1}$ be a sequence of nonnegative
   weights.
We showed under some Lyapunov-type mean-reverting assumptions  on the
coefficients  of the SDE and some light conditions on the steps and on the  weights that,
\begin{equation}\label{convarticle1}
\bar{\nu}_n(\omega,f)=\frac{1}{\eta_1+\ldots+\eta_n}\sum_{k=1}^n \eta_k f(\bar{X}_{k-1}(\omega))\xrightarrow{n\rightarrow+\infty}\nu(f)\qquad a.s.,
\end{equation}
for a large class of functions $f$ including bounded continuous
functions (see Proposition \ref{principal} below, or \cite{panloup} for more general results, $e.g.$,
when $\nu$ is not unique).  We obtained this result for two types of Euler schemes:
the \guil exact" Euler scheme that is built using the true increment of the Lévy process and some
\guil approximate" Euler schemes in
which the Lévy process increments
are replaced by an   approximation which can be simulated.\\
The aim of this paper is to study the  rate of $a.s.$ weak
convergence of  $(\bar{\nu}_n)$ toward $\nu$ for these schemes and
to devise some variants of our schemes which speed up this rate. This problem has  been first studied,
for strongly mean-reverting Brownian diffusions, by
Lamberton and Pagès (\cite{bib2}) when $\eta_n=\gamma_n$, and
by Lemaire (\cite{lemaire}) for more general weight sequences (see also
\cite{bib3} and \cite{bib25}). In particular, Lemaire
established in \cite{lemaire}  that considering some more general weights does not improve
the  rate
obtained with $\eta_n=\gamma_n$ (although some choices may improve the \guil sharp" rate). Following this remark and in order
to limit the technical difficulties, we will focus on the
case $\eta_n=\gamma_n$. However, we will assume that $(X_t)$ solution to the
Lévy driven SDE (see \eqref{edss}) is a weakly mean-reverting stochastic process,
$i.e.$ that $(X_t)$ satisfies a
weaker Lyapunov assumption than in the previously cited papers. \\
\noindent
As a first result, we show that the rate induced by the Exact
Euler scheme (Scheme (E)) is the same as that obtained for Brownian diffusions, provided
the Lévy process
 has moments up to order 4. In particular, the best rate
 is of order $n^\frac{1}{3}$ (see Theorems  \ref{vitesseAbis} and \ref{vitesseA}).
 However, in practice,  this \guil exact" scheme needs
  the increments
 of the jump component of the Lévy process to be simulated in an exact way. This  is not possible
in general  except in some particular cases
(stable processes, compound Poisson process, Gamma processes,\ldots).
 That is why we need to  consider some
approximate  Euler schemes
built with some approximations of the jump component, especially when the Lévy process jumps
infinitely often on any compact time interval.\\
The  canonical way to approximate the jump component is to
truncate its small jumps (Scheme (P)). This amounts to replacing this jump
component by a compensated compound Poisson process (CCPP). For
this type of approximation, the smaller the truncation threshold is,
the closer  the law of the corresponding CCPP is  to that of
the true jump component, but conversely, the higher the intensity of its jumps is.
So, there is a  conflict between the approximation of the jump component increments and the complexity of its
simulation procedure (when there are too many jumps). The choice of the truncation threshold is the result
of a compromise between these  constraints.
It is time varying depending on  the sequence $(\gamma_n)$ and on the Lévy measure.
When the jump component has
integrable variation, we show that it is possible to find a compromise which preserves the best rate
of the exact Euler scheme.
We mean that it is possible to construct a step sequence
$(\gamma_n)$ and a sequence of truncation thresholds such that on the one hand,
the best rate induced by this type of approximation is of order $n^\frac{1}{3}$ (see Proposition \ref{vitpart})
and on the other hand the mean  number of jumps at each time step remains uniformly bounded.
This implies that the algorithm has a linear
 mean-complexity. Otherwise, this constraint of simulation slows down
 the best achievable rate. In particular,
 when the local behavior of the jump component is very irregular, Scheme (P)
 provides some very slow rates of convergence.\\
 We propose to overcome this problem by adapting a work by Asmussen and Rosinski (\cite{bib10}) in which it is
 shown that when the truncation threshold tends to 0, the small jump component of
 a one-dimensional Lévy process has asymptotically a Brownian
 behavior. It  can be extended to $d$-dimensional
Lévy processes (see Cohen and Rosinski, \cite{cohen}).
We then construct another Euler scheme (see Scheme (W)) by
a \textit{wienerization of the small jumps}. For this scheme, the compromise
between the simulation and the approximation of the jump component is less constraining. Actually,
we show that if the
jump component has 3/2-integrable variation, it is possible to preserve the  rate of $n^\frac{1}{3}$
and to respect the constraint of simulation. Furthermore, if $\pi$
is symmetric in a neighborhood of 0, the preceding assertion is
valid without any conditions
on the small jumps (see Theorem \ref{vitesseWbis} and Proposition \ref{vitpart}).

\noindent
Before outlining the structure of the paper, we list some notations:\\
\noindent
$\bullet$ The set  $\mathbb{M}_{d,l}$ of matrices with $d$ rows and $l$ columns and real-valued entries will be
endowed with the norm  $\|M\|:=
\sup_{\{|x|\le 1\}}|Mx|/|x|$.\\
\noindent
$\bullet$ For $x\in \ER^d$ and $k\in\EN$, $x^{\otimes k}$ denotes the element of $(\ER^d)^k$ defined by
 $x^{\otimes k}_{i_1,\ldots,i_k}= x_{i_1}x_{i_2}\ldots x_{i_n}$ for every $i_1,\ldots i_k\in\{1,\ldots,k\}$.\\
 \noindent
$\bullet$ For every ${\cal C}^k$-function $f:\ER^d\mapsto\ER$ and $x, y\in\ER^d$, we  adopt the following notation :
$$D^k f(x)y^{\otimes k}=\sum_{i_1,\ldots,i_k\in\{1,\ldots,k\}}\frac{\partial^kf}{\partial x_{i_1}\ldots
\partial x_{i_k}}(x)y_{i_1}\ldots y_{i_k}.$$
If $D^k f$ is bounded, we set
$$\|D^k f\|_\infty=\sup_{i_1,\ldots,i_k\in\{1,\ldots,k\}}
\sup_{x\in\ER^d}|\frac{\partial^kf}{\partial x_{i_1}\ldots
\partial x_{i_k}}(x)|.$$
\noindent $\bullet$ We say that $V:\ER^d\mapsto \ER_+^*$ is an  EQ-function
(for \textit{Essentially Quadratic} function)  if $V$ is  a ${\cal
C}^2$-function such that $\lim V(x)=+\infty$ when $|x|\rightarrow+\infty$,
$\abs[\nabla V]\le C \sqrt{V}$ and ${D^2}V$
 is bounded.\\
\noindent $\bullet$ We set $\Gamma_n=\sum_{k=1}^n\gamma_k$, and for $s>0$, $\Gamma_n^{(s)}=\sum_{k=1}^n\gamma_k^{(s)}$.\\

\noindent In Section \ref{section1}, we introduce the framework
and the algorithm, and we recall a result of convergence
of the sequence of empirical measures established in \cite{panloup}.
In Section \ref{section2}, we state our main results
about the rate of convergence induced by  the exact and  approximate Euler
schemes when the Lévy process has moments higher than 4. Sections \ref{proofvitesseA} (resp.
\ref{proofvitesseW}) are devoted to the proof of these results in
the exact case (resp. approximate case). In Section \ref{extension}, we state a partial  extension
of the main results when the Lévy process has less moments. Finally, in Section
\ref{simulations}, we propose some numerical illustrations of  our theoretical results.
\section{\large{Setting and Background on convergence results}}\label{section1}
For a Lévy measure $\pi$ on $\ER^{l}$, we denote by $\bf{(H_p)}$ the following moment assumption
\begin{equation*}
{\bf{(H_p)}}\quad:\quad\int_{|y|>1}|y|^{2p}\pi(dy)<+\infty\quad\textnormal{with }p\ge 1.
\end{equation*}
We recall that a Lévy process with Lévy measure $\pi$ is   $2p$-integrable (see $e.g.$ \cite{bib18}, Theorem 6.1).
In \cite{panloup}, we  studied the convergence to the invariant measure for every $p>0$.
Here, we only consider the $p\ge 1$ case
because our main problem is to observe the impact of the approximation of the jump component which only depends
on the small jumps.\\
\noindent
Throughout this paper, we denote by $(X_t)_{t\ge 0}$ a solution to the following SDE
\begin{equation}\label{edss}
dX_{t}=b(X_{t^{-}})dt+\sigma(X_{t^{-}})dW_{t}+\kappa(X_{t^{-}})dZ_{t}
\end{equation}
where $b:\ER^d\mapsto\ER^d$, $\sigma:\ER^d\mapsto\mathbb
{M}_{d,l}$  and
$\kappa:\ER^d\mapsto\mathbb {M}_{d,l}$ are continuous  with
sublinear growth, $(W_{t})_{t\ge0}$ is a   $l$-dimensional Brownian motion and
$(Z_{t})_{t\ge0}$ is a    locally square-integrable purely discontinuous $\ER^{l}$-valued Lévy
process independent of  $(W_{t})_{t\ge0}$ with Lévy measure $\pi$ and characteristic
function given for every $t\ge 0$ by
$$\ES\{e^{i<u,Z_t>}\}=\exp\big[t\big(\int e^{i<u,y>}-1-i<u,y>\pi(dy)\big)\big].$$
We recall that  $(Z_t)_{t\ge 0}$ is a CCPP if and only if $\pi$ is a finite
measure and that, otherwise, it can be constructed as a limit of CCPP:
let $(u_n)_{n\ge1}$ be a sequence of positive numbers converging to 0. Let $D_n=\{|y|>u_n\}$ and
let $((Z_{t,n})_{t\ge 0})_{n\ge 1}$ denote the sequence of
processes defined by
\begin{equation}\label{sautapprox}
Z_{t,n}:=\sum_{0<s\le t}\Delta Z_s 1_{\{\Delta
  Z_s \in D_n\}}-t\int_{D_n}y\pi(dy)\qquad\forall t\ge0.
\end{equation}
For every $n\ge1$, $(Z_{t,n})_{t\ge 0}$ is a CCPP with intensity $\lambda_n=\pi(D_n)$
and jump size distribution $\mu_n(dx)=1_{D_n}\frac{\pi(dx)}{\pi(D_n)}$.
Furthermore, $Z_{.,n}\xrightarrow{n\rightarrow+\infty}Z$  in $L^2$
locally uniformly, $i.e.$
$$\ES\big\{\sup_{0\le t\le
T}|Z_t-Z_{t,n}|^2\big\}\xrightarrow{n\rightarrow+\infty}0\qquad
\forall T>0.$$ \noindent \textbf{Discretization of the SDE.} We
introduce three Euler schemes. Scheme (E)
 is constructed with the exact increments of the jump component and is called  the exact Euler scheme.
 Schemes (P) and (W) are approximate Euler schemes. In scheme (P), we truncate the small jumps
 and in scheme (W), we refine the approximation by a \textit{wienerization} of the small jumps.\\
\noindent
Let $(\gamma_n)_{n\ge1}$ be a decreasing sequence  of
 positive numbers such that
${\lim} \gamma_n=0$ and such that $\Gamma_n\rightarrow+\infty$.
Let $(U_n)_{n\ge1}$ be a sequence of i.i.d. square
integrable centered $\ER^l$-valued random variables  such
that  $\Sigma_{U_1}=I_l$. Finally, let
$(\bar{Z}_n)_{n\ge1}$, $(\bar{Z}^{^P}_n)_{n\ge1}$ and
$(\bar{Z}^{^W}_n)_{n\ge1}$ be  sequences  of independent
$\ER^l$-valued random variables, independent of $(U_n)_{n\ge1}$
satisfying
\begin{equation*}
\bar{Z}_n\overset{\cal L}{=}Z_{\gamma_n},\qquad \qquad
\bar{Z}^{^P}_n\overset{\cal
L}{=}Z_{\gamma_{n},n}\quad\textnormal{and}
\quad\bar{Z}^{^W}_n\overset{\cal
L}{=}\bar{Z}_n^{^P}+\sqrt{\gamma_{n}}Q_n\Lambda_{n}\qquad \forall
n\ge1,
\end{equation*}
where $(\Lambda_n)_{n\ge 1}$ is a sequence of i.i.d. random
variables, independent of $(\bar{Z}^{^P}_n)_{n\ge 1},U_n)_{n\ge1}$, such that  $\ES \Lambda_1=0$,
$\Sigma_{\Lambda_1}=I_d$ and $\ES\{\Lambda_1^{\otimes3}\}=0$, and
$(Q_n)$ is a sequence of $l\times l$ matrices such that
\begin{equation*}
(Q_nQ_n^*)_{i,j}=\int_{|y|\le u_k} y_i y_j\pi(dy).
\end{equation*}
We then denote by $(\bar{X}_{n})$, $(\bar{X}^{^{P}}_{n})$ and $(\bar{X}^{^{W}}_{n})$, the Euler schemes
recursively defined  by $\bar{X}_0=\bar{X}^{^P}_0=\bar{X}^{^W}_0=x\in\ER^d$ and
\begin{align*}
&{\bar{X}}_{n+1}=\bar{X}^{}_{n}+\gamma_{n+1}
b(\bar{X}^{}_{n})+\sqrt{\gamma_{n+1}}\sigma(\bar{X}^{}_{n})U_{n+1}+\kappa(\bar{X}^{}_{n})\bar{Z}_{n+1}
&\bf{(E)}\\
&\bar{X}^{^{P}}_{n+1}=\bar{X}^{^{P}}_{n}+\gamma_{n+1}
b(\bar{X}^{^{P}}_{n})+\sqrt{\gamma_{n+1}}\sigma(\bar{X}^{^{P}}_{n})U_{n+1}+\kappa(\bar{X}^{^{P}}_{n})
\bar{Z}^{^{P}}_{n+1}
&\bf{(P)}\\
&\bar{X}^{^{W}}_{n+1}=\bar{X}^{^W}_{n}+\gamma_{n+1}
b(\bar{X}^{^W}_{n})+\sqrt{\gamma_{n+1}}\sigma(\bar{X}^{^W}_{n})U_{n+1}+
\kappa(\bar{X}^{^W}_{n})\bar{Z}^{^W}_{n+1}.&\bf{(W)}
\end{align*}
We  denote by $({\cal F}_n)$, $({\cal F}^{^P}_n)$ and $({\cal F}^{^W}_n)$ the natural filtrations
induced by $(\bar{X}_n)$, $(\bar{X}^{^P}_n)$ and $(\bar{X}^{^W}_n)$ respectively.
\begin{Remarque} Note that $\bar{Z}_n^{^P}$ can be simulated if both the
intensity and the jump distribution of
$(Z_{t,n})_{t\ge 0}$ can be computed. Its simulation time
depends on the number of jumps of $(Z_{t,n})_t$ on $[0,\gamma_n]$.
Its mean is $\pi(D_n)\gamma_n$.
In order to ensure the linear mean-complexity of the algorithm, we  ask in practice these means
to be bounded, $i.e.$
\begin{equation}\label{condsim}
\sup_{n\ge 1}\pi(D_n)\gamma_n<+\infty.
\end{equation}
\noindent In scheme (W), $Q_n$ can be computed by the Choleski
method as an upper triangular matrix if $Q_nQ_n^*$ is definite.
 Otherwise, we can compute the principal square root of
$Q_nQ_n^*$.
\end{Remarque}
\noindent The associated sequences of empirical measures are
defined by
\begin{equation}
\bar{\nu}^{}_n=\frac{1}{H_n}\sum_{k=1}^{n}\eta_k\delta_{\bar{X}_{k-1}}\quad
\bar{\nu}^{^P}_n=\frac{1}{H_n}\sum_{k=1}^{n}\eta_k\delta_{\bar{X}^{^P}_{k-1}}
\quad\textnormal{and}\quad\bar{\nu}^{^W}_n=\frac{1}{H_n}\sum_{k=1}^{n}\eta_k\delta_{\bar{X}^{^W}_{k-1}}
\end{equation}
where  $(\eta_k)$ is a sequence of positive numbers such that
$H_n=\sum_{k=1}^{n}\eta_k\xrightarrow{n\rightarrow+\infty}+\infty$.
\begin{Remarque} As already mentioned, the rate of convergence will be only studied
in the case $\eta_k=\gamma_k$. However, in
the proof, we will intensively make use of  convergence results for more general weighted
empirical measures. That is why
Proposition
\ref{principal} is recalled in quite a general setting.
\end{Remarque}
\noindent Let  us pass now to the Lyapunov mean-reverting assumption.
Let $a\in
(0,1]$ be a parameter relative to  the
mean-reversion intensity. Let $r\ge0$  be a parameter relative to  the growth of the noise coefficients $\sigma$
and $\kappa$. In the sequel, we assume that there exists an $EQ$-function $V$ such that
\begin{align*}&\textbf{Assumption }
\mathbf{(S_{a,r}):}& |b|^2\le CV^a \quad
{\rm{Tr}}(\sig)+\|\kappa\|^{2}\le C V^{r}\;\textnormal{ with $r<a$.}\\
&\textbf{Assumption } \mathbf{(R_{a}):}&
\psg\nabla V,b\psd\le \bar{\beta}- \bar{\alpha} V^a \quad \textnormal{
with $\bar{\alpha}>0$ and $\bar{\beta}\in\ER$.}\quad
\end{align*}
The first deals with the growth control of the coefficients and the second is called
the mean-reverting assumption. These two assumptions  imply
assumptions $\mathbf{(S_{a,p,q})}$ and $\mathbf{(R_{a,p,q})}$
introduced in \cite{panloup}. Hence, we derive the following result from \cite{panloup}:
\begin{prop}
\label{principal} Let  $a\in(0,1]$, $p\ge 1$ and $r\in[0,a)$. Assume
 $\mathbf{(H_p)}$, $\mathbf{(R_{a})}$ and
$\mathbf{(S_{a,r})}$. Assume
$\ES\{|U_1|^{2p}\}+\ES\{|\Lambda_1|^{2p}\}<+\infty$ and
$(\eta_n/\gamma_n)$ nonincreasing. \\
\noindent
(a) i. Then,
\begin{equation}\label{Vsuptendu}
\sup_{n\ge 1}\bar{\nu}^{}_n(V^{\frac{p}{2}+a-1})<+\infty\quad a.s.
\end{equation}
Hence, the sequence $(\bar{\nu}^{}_n)_{n\ge
1}$ is $a.s.$ tight as soon as $p/2+a-1>0$.\\
ii. Moreover, if
$\kappa(x)\overset{|x|\rightarrow+\infty}{=}o(|x|)$ and
${\rm{Tr}}(\sig)+\|\kappa\|^{2}\le C V^{\frac{p}{2}+a-1}$, then every
weak limit  of $(\bar{\nu}_n)$ is an invariant probability for the
SDE \eqref{edss}. In particular,  if $(X_t)_{t\ge 0}$ admits a
unique invariant probability $\nu$, then for every continuous function
$f$ such that $f=o(V^{\frac{p}{2}+a-1})$,
$\underset{n\rightarrow\infty}{\textnormal{lim}}\bar{\nu}^{}_n(f)=\nu(f)$.\\
iii. Furthermore, $\ES\{V^p(\bar{X}^{}_n)\}=O(\Gamma_n)$ and if $a=1$, $\sup_{n\ge
1}\ES\{V^p(\bar{X}^{}_n)\}<+\infty$.\\
\noindent (b) The same result holds  for
$(\bar{\nu}^{^P}_n)_{n\ge1}$ and $(\bar{\nu}^{^W}_n)_{n\ge1}$.
\end{prop}
\begin{Remarque} For schemes (E) and (P), the above proposition is a direct consequence of Theorem 2
and Proposition 2 of \cite{panloup}.
 We did not study scheme (W) in \cite{panloup} but it is straightforward to show that the
 proposition holds true with a similar proof as that used for scheme (P).
 \noindent
Note that when $a=1$, $n\mapsto\ES\{V^p(\bar{X}_n)\}$ is bounded whereas when $a<1$, $i.e.$ when the intensity
of the mean-reverting is weak, one only has a control of its growth.
This induces some technicalities
 but has no significant influence on the main results.
 \end{Remarque}
\section{\large{Main results}} \label{section2} In this section, we suppose that
$\ES|Z_t|^{2p}<+\infty$ with $p>2$ (see Section \ref{extending} for an extension
to $p\in[1,2)$). Let $A$ denote the infinitesimal generator of $(X_t)$. $A$ is given
for every ${\cal C}^2$-function $f$ with bounded second derivatives\footnote{
Note that for such function, $Af$ is well-defined since $\ES|Z_t|^2<+\infty$.} by,
\begin{equation*}
\begin{split}
Af(x)&=\psg\nabla f,b\psd(x)+\frac{1}{2}Tr(\sigma^* D^2 f\sigma)(x)\\
& +\int \left(f(x+\kappa(x)y)-f(x)-\psg\nabla
f(x),\kappa(x)y\psd1_{\{|y|\le 1\}}\right)\pi(dy). \end{split}
\end{equation*}
We evaluate the rate of convergence
on some  test functions $g$ such that  $g=Af+C$
where  $C$ is a nonnegative real number and
 $f$ satisfies the following assumption:\\

 \noindent
 $\mathbf{(C^p_f) :}$ $(i)$  $f\in {\cal C}^4(\ER^d)$ and   $f(x)= O(V(x))$
 as $|x|\rightarrow+\infty$.\\
\noindent
$(ii)$ For $k=2,3,4,$ $D^k f$ is a bounded and Lipschitz function.\\
\noindent
$(iii)$ $|\nabla f(x)|^2=O(V^\frac{\epsilon}{2}(x))$ as $|x|\rightarrow+\infty$ with
$\epsilon\in [0,p/2 +a-1-r)$.\\

\noindent
Since $\nu$ is invariant for the SDE \eqref{edss}, we know that
$\nu(Af)=0$ (see $e.g.$ \cite{bib14})
 and then, $\nu(g)=\nu(Af+C)=C$. It follows that it suffices to evaluate
 the rate when $C=0$.\\
 \noindent
 \begin{Remarque} For a jump diffusion like \eqref{edss}, we are not able to characterize
  simply the set of functions $g$
 which can be represented as  $g=Af+C$ with $f$ satisfying
 $\mathbf{(C_f^p)}$. However, in the case of Brownian diffusions processes,
   some important works  have been done in that direction.
 Actually, in \cite{Parver1}, \cite{Parver2} and \cite{Parver3}, Pardoux and Veretennikov show that in a Sobolev framework, existence
 and unicity hold for the Poisson equation $g-\nu(g)=Af$ where $A$ is
 the infinitesimal generator of a positive recurrent diffusion. Moreover, in \cite{bib2}, Lamberton and Pagès show
 that when the diffusion is an Ornstein-Uhlenbeck process, the above equation can be solved
  in ${\cal C}^2(\ER^d)$.
 \end{Remarque}
For this class of functions, the global structure of the rates of convergence  is elucidated. For Scheme (E), our main result
is Theorem \ref{vitesseAbis}. We show that for every sequence $(\gamma_n)$, there exists a sequence
$(\rho_n)$ such that $(\rho_n\bar{\nu}_n^{^{E}}(Af))$  converges weakly:
a fast-decreasing sequence $(\gamma_n)$ (in a sense being precised in Theorem \ref{vitesseAbis}$(a)$)
leads to a CLT and a slowly-decreasing sequence $(\gamma_n)$  leads to a convergence in probability to
a deterministic constant (see Theorem \ref{vitesseAbis}$(b)$).  The rate $(\rho_n)$ is maximal for a \guil critical"
choice
of $(\gamma_n)$ for which
both types of convergence occur simultaneously. In particular, if $\gamma_n=\gamma_1 n^{-\zeta}$ with
$\zeta\in(0,1]$,
the  best rate holds
for $\zeta={\frac{1}{3}}$ (see \guil Particular Case"). In this case, $\rho_n$ is of order $n^{\frac{1}{3}}$.\\
As concerns the approximate Euler schemes,  our main results are Theorem \ref{vitesseWbis} and Proposition \ref{vitpart}.
 In the first one, we describe the structure of the rate induced by Scheme (P) and (W) as a function of
 $(\gamma_n)$ and  of $(u_n)$. When $(u_n)$ decreases \guil sufficiently fast" in a sense depending
 on the choice of the scheme, on $(\gamma_n)$ and on the Lévy measure,
 the  result induced by Scheme (E) remains valid for schemes (P) and (W). Otherwise, the approximation of
 the jump component
 dictates
 a slower rate of convergence. \\
 Theorem \ref{vitesseWbis} can not be directly applied in practice because it does not
 specify whether the fundamental condition of simulation \eqref{condsim} is compatible
 with the theoretical results.
 This is the purpose of Proposition \ref{vitpart} in which we give the best possible
 rates for schemes (P) and (W)
 under condition \eqref{condsim} as a function depending on the local behavior of the small jumps.
 In particular, Proposition  \ref{vitpart} clarifies   the impact  of the
wienerization
of the small jumps announced in the introduction and shows that it makes possible to preserve
 the same rate of convergence of the exact Euler scheme for a wide class of Lévy processes (for
 which the exact simulation of the increments is impossible).\\

\begin{figure}[h]
\centering
\includegraphics[width=7cm]{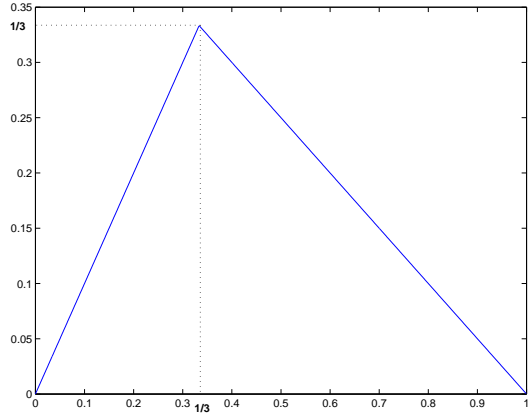}
\caption{Rate of convergence for polynomial steps} \label{tente}
\end{figure}

\noindent Let $f\in C^1(\ER^d)$. We define $\tilde{H}^f$ by
$$\tilde{H}^f(z,x,y)=f(z+\kappa(x)y)-f(z)-\psg\nabla f(z),\kappa(x)y\psd$$
and $z\mapsto\tilde{H}^f(z,x,y)$ is denoted by  $\tilde{H}^f_{.,x,y}$. Our first result is the following:
\begin{theorem}\label{vitesseAbis} Assume that $\ES|Z_t|^{2p}<+\infty$ with $p>2$ and that \eqref{edss} admits a
unique
 invariant measure $\nu$.
Let $a\in(0,1]$ and $r\ge 0$ such  that
 $\mathbf{(R_{a})}$  and
$\mathbf{(S_{a,r})}$ are satisfied and  $p/2+a-1> 2r$.  If moreover, $\ES\{U_1^{\otimes3}\}=0$,
$\ES\{|U_1|^{2p}\}<+\infty$ and
$\eta_n=\gamma_n$ for every $n\ge 1$, then for every function $f:\ER^d\mapsto\ER$
satisfying $\mathbf{(C^p_f)}$,

\noindent(a) If
$\frac{\Gamma_n^{(2)}}{\sqrt{\Gamma_n}}\xrightarrow{n\rightarrow+\infty}\hat{\gamma}\in
[0,+\infty)$, $\sqrt{\Gamma_n}\bar{\nu}^{}_n(Af)\underset{n\rightarrow+\infty}{\xrightarrow{\cal
L}}{\cal N}\Big(\hat{\gamma}m,\hat{\sigma}_f^2 \Big).$\\
(b) If
$\frac{\Gamma_n^{(2)}}{\sqrt{\Gamma_n}}\xrightarrow{n\rightarrow+\infty}+\infty$,
$\frac{\Gamma_n}{\Gamma_n^{(2)}}
\bar{\nu}^{}_n(Af)\underset{n\rightarrow+\infty}{\overset{\PE}{\longrightarrow}}m$\\
\noindent
where $\hat{\sigma}_f^2=\int\big(|\sigma^*\nabla f|^2(x)+\int
\big(f(x+\kappa(x)y)-f(x)\big)^2\pi(dy)\big)\nu(dx)$  and,
\begin{align*}
&m=-\int(\phi_1(x)+\phi_2(x)+\phi_3(x))\nu(dx) \textnormal{ with }\phi_1(x)=\frac  1 2  {D^2}f(x)b(x)^{\otimes2},\\
&\phi_2(x)= \int\frac 1 6 D^3 f(x); b(x) ;(\sigma(x)u)^{\otimes 2}+\frac{1}{24}
D^4 f(x)(\sigma(x)u)^{\otimes 4}\PE_{U_1}(du)\\
\textit{and,}\quad&\phi_3(x)=\frac{1}{2}\int\pi(dy_1)\int\pi(dy_2)\tilde{H}^{\tilde{H}^f_{.,x,y_1}}(x,x,y_2)\\
&+\int\pi(dy_1)\Big(\psg\nabla\tilde{H}^f_{.,x,y_1}(x),b(x)\psd+\int\PE_{U_1}(du){D^2}(\tilde{H}^f_{.,x,y_1})
(x)(\sigma(x)u)^{\otimes2}\Big).
\end{align*}
\end{theorem}
\noindent
\textbf{Particular Case.} Assume that $\gamma_n=\gamma_1n^{-\zeta}$ with $\zeta\in(0,1]$. Then,
\begin{eqnarray*}
\begin{cases}
\sqrt{\gamma_1\log n}\bar{\nu}^{}_n(Af)\underset{n\rightarrow+\infty}{\xrightarrow{\cal
L}}{\cal N}\Big(0,\hat{\sigma}_f^2 \Big).&\textnormal{if $\zeta=1$}\\
\sqrt{\frac{\gamma_1}{1-\zeta}}n^{\frac{1-\zeta}{2}}\bar{\nu}^{}_n(Af)\underset{n\rightarrow+\infty}{\xrightarrow{\cal
L}}{\cal N}\Big(\hat{\gamma}m,\hat{\sigma}_f^2 \Big).&\textnormal{if $\zeta\in[1/3,1)$}\\
\frac{1-2\zeta}{\gamma_1(1-\zeta)}n^\zeta
\bar{\nu}^{}_n(Af)\underset{n\rightarrow+\infty}{\overset{\PE}{\longrightarrow}}m&\textnormal{if $\zeta<1/3$}
\end{cases}
\end{eqnarray*}
where $\hat{\gamma}=0$ if $\zeta\in(1/3,1)$ and $\hat{\gamma}=\sqrt{6\gamma_1^3}$ if $\zeta=1/3.$
On Figure \ref{tente}, one represents $\zeta\mapsto h(\zeta)$
 where $h(\zeta)$ denotes the exponent of the rate. One observes that $\underset{\zeta\in(0,1]}{\max} h(\zeta)=h(1/3)=1/3$.
\begin{Remarque} Theorem \ref{vitesseAbis} shows that the rate is  the same as that obtained
for Brownian diffusions. In particular, when $\kappa=0$, Theorem \ref{vitesseAbis}
extends the rate  results of \cite{bib2} and \cite{lemaire}  to the
weakly mean-reverting diffusions ($a<1$), whose convergence to the invariant measure
has been studied in \cite{bib3}.\\
\noindent
Note that the condition $\ES\{U_1^{\otimes3}\}=0$ is not necessary for the convergence of the empirical
measures but plays a role in the rate. Without this condition, the best rate would be
of order $n^{\frac{1}{4}}$, obtained for $\zeta=1/2$ (see \cite{bib2} in the case of Brownian diffusions).
\end{Remarque}
\noindent
Let us pass now to the main results for the approximate Euler schemes. Let $(u_k)_{k\ge 1}$ denote the sequence of
truncation thresholds and set $ {\cal
\beta}_{n,\pi}^{(s)}=\sum_{k=1}^n \gamma_k \int_{|y|\le
u_k}|y|^s\pi(dy)$.
For $s\in\{2,3,4\}$, we introduce a new assumption ${\bf{(A_s^1)}}$ which is relative to the  impact of
the jump component approximation
as a function of the steps and of  the truncation
thresholds:
\begin{equation*}
{\bf{(A_s^1):}}\quad\frac{\beta_{n,\pi}^{(s)}}{\Gamma_n^{(2)}}\xrightarrow{n\rightarrow+\infty}\hat{\alpha}_s\in[0,+\infty]
\quad
\textnormal{and}\quad
\frac{\beta_{n,\pi}^{(s)}}{\sqrt{\Gamma_n}}\xrightarrow{n\rightarrow+\infty}\hat{\beta}_s\in[0,+\infty].
\end{equation*}
Since $s\mapsto\beta_{n,\pi}^{(s)}$ is a decreasing function,
 $s\mapsto\hat{\alpha}_s$ and
 $s\mapsto\hat{\beta}_s$ both decrease. This can be interpreted as follows: the constraint
 on $(u_k)$ decreases with $s$. \\
For $s\in\{2,3,4\}$,
we also introduce another assumption on the Lévy measure that will be necessary to transform some tightness results
 in some  convergence in distribution results:
 \begin{equation*}
 {\bf{(A_s^2)}}:\textnormal{ For every $i_1,\ldots,i_s\in\{1,\ldots,l\}$, }
\Big(\frac{\int_{|y|\le u_k} y_{i_1}\ldots y_{i_s}\pi(dy)}{\int_{|y|\le
u_k} |y|^{ s}\pi(dy)}\Big)_{k\ge 1} \textnormal{ converges in }\ER .
\end{equation*}
For instance, the above assumption is satisfied if $\pi(dy)=\psi(y)\lambda_l(dy)$ where $\psi$ satisfies:
there exists $\alpha\in[l,l+2)$ such that $|y|^{\alpha+l}\psi(y)\rightarrow C_0\in\ER_+^*$ when $y\rightarrow0$. \\
\noindent
Throughout this paper, we will say that the Lévy measure \textit{$\pi$ is  quasi-symmetric in a neighborhood of 0} if
  $\int_{\{|y|\le u\}} y^{\otimes3}\pi(dy)=0$ for $u$
sufficiently small. In particular, this assertion holds if $\pi$ is symmetric in a neighborhood of 0. \\
\noindent
We will also say that a real-valued random variable $X$ is \textit{quasi-subgaussian} if there exists
$m>0$ and $\sigma>0$ such that for all
$M>0$
$$ \PE(|X|>M)\le \PE(|Y|+m>M)\quad \textnormal{with}\quad Y\sim{\cal N}(0,\sigma^2).$$
\begin{theorem} \label{vitesseWbis}
Let $a\in(0,1]$, $r\ge 0$, $p<2$ such that the conditions of Theorem \ref{vitesseAbis}
are satisfied. Assume that $\ES\{|\Lambda_1|^{2p}\}<+\infty$ and
that $(u_k)_{k\ge1}$ is decreases to 0.\\
(a) i. Scheme (P): Assume that $\mathbf{(A_s^1)}$ holds with $s=2$.\\
\noindent
$\bullet$ If $\hat{\alpha}_s=0$ or $\hat{\beta}_s=0$, then the conclusions of Theorem \ref{vitesseAbis} are
still valid for $(\bar{\nu}_n^{^P})$.\\
$\bullet$ If $\hat{\alpha}_s\in(0,+\infty]$ and $\hat{\beta}_s\in(0,+\infty)$, then
$(\frac{\Gamma_n}{\beta_{n,\pi}^{(s)}}\bar{\nu}^{^P}_n(Af))_{n\ge
 1}$ is tight with quasi-subgaussian limiting distributions.\\
$\bullet$ If $\hat{\alpha}_s\in(0,+\infty]$ and $\hat{\beta}_s=+\infty$, then
$(\frac{\Gamma_n}{\beta_{n,\pi}^{(s)}}\bar{\nu}^{^P}_n(Af))_{n\ge
 1}$ is tight with bounded limiting distributions.\\
 \noindent
 ii. Scheme (W): Assume that $\mathbf{(A_3^1)}$ holds.\\
\noindent Then, the conclusions of $(a).i$ are valid for $(\bar{\nu}^{^W}_n(Af))_{n\ge
 1}$ with $s=3$.
\noindent
Furthermore, if $\pi$ is quasi-symmetric in a neighborhood of 0 and $\mathbf{(A_4^1)}$ holds,
the conclusions of $(a).i$ are valid for $(\bar{\nu}^{^W}_n(Af))_{n\ge
 1}$ with $s=4$.\\
 \\
\noindent
(b) i. Scheme (P): Assume that $\mathbf{(A_s^1)}$ and $\mathbf{(A_s^2)}$ hold with $s=2$. Then,
\begin{align*}
&\frac{\Gamma_n}{\beta_{n,\pi}^{(s)}}\bar{\nu}^{^P}_n(Af)
 \underset{n\rightarrow+\infty}{\xrightarrow{\cal
L}}{\cal N}\Big(m/\hat{\alpha}_s -m_s,(\hat{\sigma}_f/\hat{\beta}_s)^2 \Big)&
\textnormal{if $\hat{\alpha}_s\in(0,+\infty]$
 and $\hat{\beta}_s\in(0,+\infty)$}\\
 &\frac{\Gamma_n}{\beta_{n,\pi}^{(s)}}\bar{\nu}^{^P}_n(Af)
 \underset{n\rightarrow+\infty}{\xrightarrow{\PE}}m/\hat{\alpha}_s-m_s
&\textnormal{if $\hat{\alpha}_s\in(0,+\infty]$
 and $\hat{\beta}_s\in(0,+\infty)$,}
 \end{align*}
with  $|m_2|\le \bar{m}_2=(d/2)\|{D^2}f\|_{\infty}\int\|\kappa\|^2(x)\nu(dx)$
and $m$ and $\hat{\sigma}_f^2$ like in Theorem \ref{vitesseAbis}.\\
ii. Scheme (W): Assume that $\mathbf{(A_3^1)}$ and $\mathbf{(A_3^2)}$ hold.\\
\noindent Then, the conclusions of $(b).i$ are valid for $(\bar{\nu}^{^W}_n(Af))_{n\ge 1}$ with $s=3$
and a real number $m_3$ satisfying $|m_3|\le \bar{m}_3=(d^{\frac 3 2}/6)\|D^3f\|_{\infty}\int\|\kappa(x)\|^2\nu(dx)$.\\
\noindent
Furthermore, if $\pi$ is quasi-symmetric in a neighborhood of 0 and if $ \mathbf{(A_4^1)}$
 and $ \mathbf{(A_4^2)}$ hold, the conclusions of $(b).i$ are valid for $(\bar{\nu}^{^W}_n(Af))_{n\ge
 1}$ with $s=4$ and a real number $m_4$ satisfying  $|m_4|\le\bar{m}_4 =(d^{2}/24)\|D^4f\|_{\infty}\int\|\kappa(x)\|^4\nu(dx)$.
\end{theorem}
\begin{Remarque}\label{extension}  Note that in the one-dimensional case,  Assumption $\bf{(A_s^2)}$ is
always satisfied   when $s=2$
 or $s=4$. In those cases,  $m_s=\frac{1}{s!}\int
f^{(s)}(x) \kappa(x)^s\nu(dx)$. If $s=3$, Assumption $\bf{(A_s^2)}$ is satisfied if
$\int_{\{|y|\le u_k}y^3\pi(dy)/\int_{\{|y|\le u_k}|y|^3\pi(dy)\rightarrow a_3\in\ER$.
In this case, $m_3=a_3\frac{1}{3!}\int
f^{(3)}(x) \kappa(x)^3\nu(dx)$. In the multidimensional case, the value of $m_s$ is also explicit
but its expression is more complicated (see proof of Lemma \ref{resteB}).
\end{Remarque}
Let us now state Proposition \ref{vitpart}. In $(a)$, we provide some conditions on the Lévy measure
in the neighborhood of 0 which preserve the rate of convergence induced by the exact Euler scheme
under the  condition of simulation \eqref{condsim}. In $(b)$, we suppose that the Lévy measure has a density
closed to that of an $\alpha$-stable process
in the neighborhood of 0 and give in that case the optimal rate for the two schemes as a function of $\alpha$.
For these two parts, we also give some available choices of steps and truncation thresholds.
\begin{prop}\label{vitpart} Let $a\in(0,1]$, $r\ge 0$, $p\ge 2$ such that the conditions
of Theorem \ref{vitesseAbis}
are satisfied. Assume that $\ES\{|\Lambda_1|^{2p}\}<+\infty$  . \\
(a) Assume  that $\int_{\{|y|\le 1\}}|y|^{q}\pi(dy)<+\infty$ with
$q\in[0,2]$ and set $\gamma_k=\gamma_1 k^{-\frac{1}{3}}$ and
$u_k=\gamma_k^r$ with
$r\in [\frac{1}{q}, \frac{1}{s-q}]$. Then, Condition \eqref{condsim} holds and, \\
i. Scheme (P): If  $q\le 1$ and $s=2$,
$n^\frac{1}{3}\bar{\nu}^{^P}_n(Af)\underset{n\rightarrow+\infty}{\xrightarrow{\cal
L}}\sqrt{2/3}{\cal N}\big(m\sqrt{6},\hat{\sigma}_f^2
\big)$.\\
\\
ii. Scheme (W): If $q\le 3/2$ and $s=3$, $n^\frac{1}{3}\bar{\nu}^{^W}_n(Af)\underset{n\rightarrow+\infty}
{\xrightarrow{\cal
L}}\sqrt{2/3}{\cal N}\big(m\sqrt{6},\hat{\sigma}_f^2
\big)$.\\
Furthermore, if $\pi$ is quasi-symmetric in the neighborhood of 0, the preceding assertion
is valid with $s=4$ and every $q\in[0,2]$.\\
(b) Assume that there exists $\epsilon_0>0$ such that
\begin{equation}\label{stablecomp}
\pi(dy)=\psi(y)\lambda_l(dy)\quad\textnormal{with} \quad
1_{\{0<|y|\le\epsilon_0\}}\frac{C_1}{|y|^{\alpha+l}}\le \psi(y)\le
\frac{C_2}{|y|^{\alpha+l}}1_{\{0<|y|\le\epsilon_0\}}.
\end{equation} Set $\gamma_k=\gamma_1 k^{-(\frac{1}{3}\vee
\frac{\alpha}{2s-\alpha})}$, $u_k=\gamma_k^r$ with $r\in
[\frac{1}{\alpha}, \frac{1}{(s-\alpha)\vee \alpha}]$. Then, Condition \eqref{condsim} holds and,\\
i. Scheme (P),
$s=2$: $\big(n^{(\frac 1 3 \wedge
\frac{2-\alpha}{4-\alpha})}\bar{\nu}^{^P}_n(Af)\big)_{n\ge1}$ is
tight.\\
ii. Scheme (W), $s=3$: $\big(n^{(\frac 1 3 \wedge
\frac{3-\alpha}{6-\alpha})}\bar{\nu}^{^W}_n(Af)\big)_{n\ge1}$ is
tight.
\end{prop}
\begin{Remarque} Figure
\ref{puisoptW} represents  $\alpha\mapsto h(\alpha)$ where
$h(\alpha)$ denotes the exponent of the optimal rate induced by each approximate scheme under
the assumptions of Proposition \ref{vitpart}$(b)$. This figure emphasizes the necessity of scheme (W) when
the jump component has infinite variation because  the optimal rate of convergence
induced by scheme (P) decreases very rapidly in that case.
\begin{figure}[h]
\centering \label{puisoptW}
\includegraphics[width=7cm]{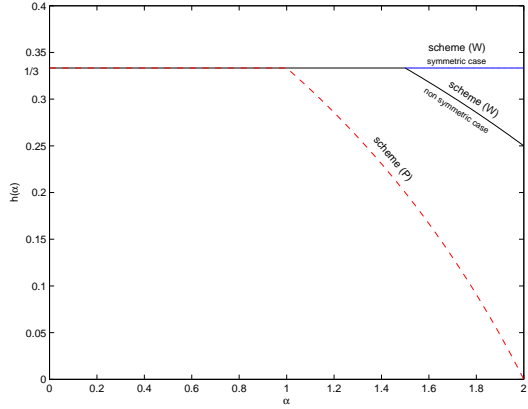}
\caption{Optimal rate in terms of the local behavior of the Lévy
process}
\end{figure}
\end{Remarque}
\begin{Remarque} The fact that we optimize the rate for the
range of sequences $(\gamma_k,u_k)_{k\ge1}$ such that the
linearity of the mean-complexity  of the procedure is ensured can
be disputable when the optimal rate is not of order $n^\frac{1}{3}$.
Actually, in this case, even if for a smaller level of truncation,
the linearity of the complexity fails, the theoretical rate is
better. Hence, another point of view consists in evaluating the
order of precision as a function of the complexity. Some precise
statements on that question would require some Berry-Esseen type
estimates (in our inhomogeneous framework). Nevertheless,
heuristic study can be done when \eqref{stablecomp} is satisfied
and suggests that the asymptotic order of precision as a function
of the mean-complexity is optimized for a class of
 steps and truncation levels including the choices of Proposition \ref{vitpart}.
\end{Remarque}
\section{\large{Proof of Theorem
\ref{vitesseAbis}}}\label{proofvitesseA} In this section, we prove the main result
induced by the exact Euler scheme: Theorem \ref{vitesseAbis}. Firstly, we decompose  $\bar{\nu}_n(Af)$
(see Lemma \ref{lemmeZ}) and then,
we compute the rate of each term of the decomposition in Lemmas
\ref{martinA}, \ref{resteA} and \ref{resteD}. (We will principally
focus  on Lemmas \ref{martinA} and \ref{resteD} where the  rate of
the jump part of the decomposition is studied). Finally, a synthesis of the previous lemmas
is realized in subsection \ref{preuvevitesseA} and completes the proof of Theorem \ref{vitesseAbis}.
\subsection{ Decomposed computation of the rate of $\bar{\nu}_n(Af)$.}\label{step1}
We set
\begin{equation*}
\bar{X}_{k,1}=\bar{X}_{k-1}+\gamma_k b(\bar{X}_{k-1}),
\quad\textnormal{and}\quad
\bar{X}_{k,2}=\bar{X}_{k,1}+\sqrt{\gamma_k}\sigma(\bar{X}_{k-1})U_k.
\end{equation*}
\noindent Denote by $(Z^{^{(k)}})_{k\ge 1}$, a sequence of
i.i.d. random variables such that $Z^{^{(1)}}\overset{\cal L}{=}Z$
and set $\bar{Z}_k=Z_{\gamma_k}^{^{(k)}}$.
\begin{lemme} \label{lemmeZ} For  $f\in { C}^2(\ER^d)$, we have
the following decomposition.
\begin{align*}
\sum_{k=1}^n\gamma_k Af(\bar{X}_{k-1})&=f(\bar{X}_n)-f(\bar{X}_0)-\sum_{k=1}^n \Big(\xi_1(\gamma_k,\bar{X}_{k-1},U_k)+\xi_2(\gamma_k,\bar{X}_{k-1},Z^{^{(k)}})\Big)\\
&-\sum_{k=1}^n
\Big(\Theta_1(\gamma_k,\bar{X}_{k-1})+\Theta_2(\gamma_k,\bar{X}_{k-1},U_k)+\Theta_3(\gamma_k,\bar{X}_{k-1},\bar{X}_{k,2},Z^{^{(k)}})\Big)\\
&-\sum_{k=1}^n
\Big((R_1+R_2)(\gamma_k,\bar{X}_{k-1},U_k)+R_3(\gamma_k,\bar{X}_{k-1},\bar{X}_{k,2},Z^{^{(k)}})\Big)\quad
\end{align*}
where,
\begin{align*}
& \xi_1(\gamma,x,U_k)=\sqrt{\gamma}\psg\nabla f(x),\sigma(x)U_k\psd ,\\
&\xi_2(\gamma,x,Z)=\int_0^\gamma\psg\nabla
f(x),\kappa(x)dZ_s\psd  +\big(\sum_{0<s\le
\gamma}\tilde{H}^f(x,x,\Delta
Z_s)-\gamma\int\tilde{H}^f(x,x,y)\pi(dy)\big),\\
&\Theta_1(\gamma,x)=\gamma\int_0^1\psg\nabla f(x+\theta\gamma
b(x))-\nabla
 f(x),b(x) \psd d\theta,\\
&\Theta_2(\gamma,x,u)=\gamma\int_0^1(1-\theta) \big({D^2}f(x+\gamma b(x)+\theta\sqrt{\gamma}\sigma(x)u)
-{D^2}f(x)\big)(\sigma(x)u)^{\otimes2}d\theta,
\end{align*}
\begin{align*}
&\Theta_3(\gamma,x,z,Z)=\sum_{0<s\le
\gamma}\big(\tilde{H}^f(x,z,Z_{s^-},\Delta
Z_s)-\tilde{H}^f(x,x,0,\Delta Z_s)\big),\\
&R_1(\gamma,x,U_k)=\sqrt{\gamma}\psg\nabla f(x+\gamma b(x))-\nabla f(x),\sigma(x)U_k\psd ,\\
&R_2(\gamma,x,U_k)=\frac {\gamma}{2} \big({D^2} f(x) (\sigma(x)U_k)^{\otimes2}-\ES\{{D^2} f(x) (\sigma(x)U_k)^{\otimes2}\}\big),\\
&R_3(\gamma,x,z,Z)=\int_0^\gamma\psg\nabla
f(z+\kappa(x)Z_{s^-})-\nabla f(x),\kappa(x)dZ_s\psd .
\end{align*}
\end{lemme}
\begin{proof}[\textnormal{\textbf{Proof}}]
We write
$$f(\bar{X}_k)-f(\bar{X}_{k-1})=\Big(f(\bar{X}_{k,1})-f(\bar{X}_{k-1})\Big)+\Big(f(\bar{X}^2_k)-f(\bar{X}^1_{k}\Big))+\Big(f(\bar{X}_k)-f(\bar{X}_{k,2})\Big)$$
We expand the first two terms by  the Taylor formula and use
the Itô formula  (with jumps) for the last one. The lemma follows
by summing up the equality for $k=1,\ldots,n.$
\end{proof}\noindent
As mentioned before, we study successively the rate of convergence
of each term of the previous decomposition. We start by showing a CLT for the terms associated
with $\xi_1$ and $\xi_2$.
\begin{lemme}\label{martinA} Assume that $\mathbf{(H_p)}$ holds for $p>2$. Let
$f:\ER^d\mapsto\ER$ satisfy   $\mathbf{(C_f^p)}$.
Then, with the notations of Lemma \ref{lemmeZ}, we have\\
(a)
\begin{equation}\label{eq1martinA}
\ES\{|\xi_2(\gamma,x,Z)|^2\}=\gamma\int
\big(f(x+\kappa(x)y)-f(x)\big)^2\pi(dy)
\end{equation}
 and there exists $\delta>0$  and a locally bounded function $C$ such  {that}\mbox{ }
\begin{equation}\label{eq2martinA}
\ES\{|\xi_2(\gamma,x,Z)|^{2(1+\delta)}\}\le C(x)\gamma
\end{equation}
(b)  Moreover, if $\mathbf{(R_a)}$ and  $\mathbf{(S_{a,r})}$
hold with  $2r<p/2+a-1$ and $\ES\{|U_1|^{2p}\}<+\infty$, then,
$$\frac{1}{\sqrt{\Gamma_n}}\sum_{k=1}^n\Big(\xi_1(\gamma_k,\bar{X}_{k-1},U_k)+\xi_2(\gamma_k,\bar{X}_{k-1},Z^{^{(k)}})\Big)\underset{n\rightarrow+\infty}{\xrightarrow{\cal L}}{\cal N}\big(0, \hat{\sigma}_f^2\big),$$
with $\hat{\sigma}^2_f=\int\Big(|\sigma^*\nabla f|^2(x)+\int
\big(f(x+\kappa(x)y)-f(x)\big)^2\pi(dy)\Big)\nu(dx)$.
\end{lemme}
\begin{proof}[\textnormal{\textbf{Proof}}]$(a)$. Let $(Z_{.,n})_{n\ge1}$ be the sequence of processes
defined by \eqref{sautapprox}. We know that
$$\ES\{\sup_{\{0\le
s\le t\}}|Z_{s,n}-Z_s|^2\}\xrightarrow{n\rightarrow+\infty}0.$$ As
$Z_{.,n}$ has bounded variations, $\xi_2(\gamma,x,Z_{.,n})$ can be
written
$$\xi_2(\gamma,x,Z_{.,n})=\sum_{0<s\le\gamma}1_{\{|\Delta
  Z_s|> u_n\}}\big( f(x+\kappa(x)\Delta Z_s)-f(x)\big)
  -\gamma\int_{\{|y|>u_n\}} \big(f(x+\kappa(x)y\big)-f(x)\big)\pi(dy).$$
Since ${D^2}f$ is bounded and $\ES|Z_t|^4<+\infty$, one easily checks that
$\xi_2(\gamma,x,Z_{.,n})$ is a locally
square-integrable purely discontinuous martingale. We deduce from
the compensation formula that
\begin{equation}\label{31'}
\ES\{|\xi_2(\gamma,x,Z_{.,n})|^2\}=\gamma\int_{\{|y|>u_n\}} |f\big(x+\kappa(x)y)\big)-f(x)|^2\pi(dy).
\end{equation}
We also check that
$$\ES\{|\xi_2(\gamma,x,Z)-\xi_2(\gamma,x,Z_{.,n})|^2\}\le
C_x \int_{\{|y|\le
u_n\}}|y|^2\pi(dy)\xrightarrow{n\rightarrow+\infty}0.$$ Letting
$n\rightarrow+\infty$ in \eqref{31'} yields the first
identity.\\
\noindent Now, let us prove the inequality. $\xi_2(\gamma,x,Z)=
\psg\nabla f(x),\kappa(x)Z_{\gamma}\psd+M_\gamma$ where $M$ is a
martingale defined by
\begin{align*}
 M_\gamma&=\sum_{0<s\le \gamma}\int_0^1 \psg\nabla
f(x+\theta\kappa(x)\Delta Z_s)-\nabla f(x),\kappa(x)\Delta
Z_s\psd d\theta\\&-\gamma\int\int_0^1\psg\nabla
f(x+\theta\kappa(x)y)-\nabla f(x),\kappa(x)y\psd d\theta\pi(dy).
\end{align*}
Let $\delta\in (0,1]$ such that $4(1+\delta)\le 2p$. Since $\nabla f$ is Lipschitz continuous, we
derive from the  Burkholder-Davis-Gundy inequality that
\begin{align*}
\ES |Z_\gamma|^{2(1+\delta)}\le\ES\Big\{\big(\sum_{0<s\le\gamma}|\Delta
Z_s|^{2}\big)^{1+\delta}\Big\}\quad\textnormal{and}\quad\ES
|M_\gamma|^{2(1+\delta)}\le C(x)\ES\Big\{\big(\sum_{0<s\le\gamma}|\Delta
Z_s|^{4}\big)^{1+\delta}\Big\}.
\end{align*}
It follows that
\begin{equation*}
\ES\{|\xi_2(\gamma,x,Z)|^{2(1+\delta)}\}\le
C_1(x)\ES\Big\{\big(\sum_{0<s\le\gamma}|\Delta
Z_s|^{2}\big)^{1+\delta}\Big\} +C_2(x)
\ES\Big\{\big(\sum_{0<s\le\gamma}|\Delta
Z_s|^{4}\big)^{1+\delta}\Big\}.
\end{equation*}
Then, it suffices  to prove that
\begin{equation}\label{1e}
\ES\Big\{\big(\sum_{0<s\le\gamma}|\Delta
Z_s|^{2+\rho}\big)^{1+\delta  }\Big\}=O(\gamma)\qquad\textnormal{
for $\rho=0$ and $\rho=2$.} \end{equation}
Denote by $(\tilde{M}_s)$ the martingale defined by
$\tilde{M}_s=\sum_{0<s\le\gamma}|\Delta Z_s|^{2+\rho}-\gamma\int|y|^{2+\rho}\pi(dy)$.
By the elementary inequality
\begin{equation}\label{elemineq}
\forall u,v\in\ER\textnormal{ and }\alpha>0,\quad|u+v|^\alpha\le 2^{\alpha\vee1-1} (|u|^\alpha+|v|^\alpha).
\end{equation}
we have,
$$\ES\Big\{\big(\sum_{0<s\le\gamma}|\Delta
Z_s|^{2+\rho}\big)^{1+\delta  }\Big\}\le C\big(\ES|\tilde{M}_\gamma|^{1+\delta}+
\gamma^{1+\delta}\int|y|^{2+\rho}\pi(dy)\big).$$
By the Burkhölder-Davis-Gundy inequality,
\begin{equation*}
\ES|\tilde{M}_\gamma|^{1+\delta}\le C\ES\Big\{\big(\sum_{0<s\le\gamma}|\Delta
Z_s|^{2(2+\rho)}\big)^{\frac{1+\delta}{2}  }\Big\}.
\end{equation*}
Since $(1+\delta)/2\le 1$, it follows from \eqref{elemineq} and from the compensation formula that
\begin{equation*}
\ES|\tilde{M}_\gamma|^{1+\delta}\le C\ES\Big\{\sum_{0<s\le\gamma}|\Delta
Z_s|^{(2+\rho)(1+\delta)}  \Big\}\le C\gamma\int|y|^{(2+\rho)(1+\delta)}\pi(dy).
\end{equation*}
Since $2\le (2+\rho)(1+\delta)\le 2p$,
 $\int |y|^{(2+\rho)(1+\delta)}\pi(dy)<+\infty$.
\eqref{1e}
follows.\\

\noindent $(b)$ Let $\{(\xi_k^n),k=1,\ldots,n,n\ge 1\}$ be a sequence of triangular
arrays of square-integrable martingale increments defined by
$$\xi_k^n=\frac{1}{\sqrt{\Gamma_n}}
\Big(\xi_1(\gamma_k,\bar{X}_{k-1},U_k)+\xi_2(\gamma_k,\bar{X}_{k-1},Z^{^{(k)}})\Big).$$
Since
$\Sigma_{U_1}=I_l$, we have
$$\ES\{|\xi_1(\gamma_k,\bar{X}_{k-1},U_k) |^2/{\cal F}_{k-1}\}=\gamma_k|\sigma^*\nabla f|^2(\bar{X}_{k-1}).$$
Moreover, $\xi_1(\gamma_k,\bar{X}_{k-1},U_k)$ and
$\xi_2(\gamma_k,\bar{X}_{k-1},Z^{^{(k)}})$ are  independent
conditionally to ${\cal F}_{k-1}$ and
$$\ES\{\xi_1(\gamma_k,\bar{X}_{k-1},U_k) /{\cal F}_{k-1}\}=\ES\{\xi_2(\gamma_k,\bar{X}_{k-1},Z^{^{(k)}}) /{\cal F}_{k-1}\}=0.$$
Then, we deduce from \eqref{eq1martinA} that
\begin{align*}
\ES\{|\xi_k^n|^2/{\cal F}_{k-1}\}&=\frac{1}{\Gamma_n}\Big(\ES\{|\xi_1(\gamma_k,\bar{X}_{k-1},U_k) |^2/{\cal F}_{k-1}\}+\ES\{|\xi_2(\gamma_k,\bar{X}_{k-1},Z^{^{(k)}})|^2/{\cal F}_{k-1}\}\Big)\\
&=\frac{\gamma_k}{\Gamma_n}\Big(|\sigma^*\nabla
f|^2(\bar{X}_{k-1})+\int
\big(f(\bar{X}_{k-1}+\kappa(\bar{X}_{k-1})y)-f(\bar{X}_{k-1})\big)^2\pi(dy)\Big).
\end{align*}
Since ${D^2}f$ is bounded, we derive from Taylor's formula and from the assumptions on
$r$ and on $\nabla f$ that
\begin{align}
\int \big(f(.+\kappa(.)y)-f(.)\big)^2\pi(dy)+|\sigma^*\nabla f|^2
\le CV^{(\epsilon+r)\vee (2r)}=o(V^{\frac{p}{2}+a-1}).\label{contri}
\end{align}
Hence, Proposition \ref{principal} yields
\begin{equation}\label{1f}
\sum_{k=1}^n\ES\{|\xi_k^n|^2/{\cal
F}_{k-1}\}\xrightarrow{n\rightarrow+\infty}\int\big(|\sigma^*\nabla
f|^2+\int \big(f(.+\kappa(.)y)-f(.)\big)^2\pi(dy)\big)d\nu.
\end{equation}
Then, the lemma will follow from the central limit
theorem for arrays of square-integrable martingale increments (see
Hall and Heyde, \cite{hall}) provided  the
\textit{Lindeberg condition} is fulfilled, $i.e.$
$$R^\rho_n=\sum_{k=1}^n\ES\{|\xi_k^n|^2 1_{\{|\xi_k^n|\ge \rho\}}/{\cal F}_{k-1}\}
\xrightarrow{n\rightarrow+\infty}0\qquad a.s.\quad \forall
\rho>0.$$ Let $A\in(0,+\infty)$ and set
\begin{align*}
&R^{\rho,A}_{n,1}=\sum_{k=1}^n1_{\{|\bar{X}_{k-1}|\le A\}}\ES\{|\xi_k^n|^2 1_{\{|\xi_k^n|\ge \rho\}}/{\cal F}_{k-1}\},\\
&R^{\rho,A}_{n,2}=\sum_{k=1}^n1_{\{|\bar{X}_{k-1}|\ge A\}}\ES\{|\xi_k^n|^2 1_{\{|\xi_k^n|\ge \rho\}}/{\cal F}_{k-1}\}.
\end{align*}
We have $\ES\{|\xi_k^n|^2 1_{\{|\xi_k^n|\ge \rho\}}/{\cal
F}_{k-1}\}=F^n_A(\bar{X}_{k-1},\gamma_k)\;\textnormal{where}$
$$F^n_A(x,\gamma)=\frac{1}{\Gamma_n}\ES\{|\xi_1(\gamma,x,U_1)+\xi_2(\gamma,x,Z)|^2
1_{\{|\xi_1(\gamma,x,U_1)+\xi_2(\gamma,x,Z)|\ge\rho\sqrt{\Gamma_n}\}}\}.$$
Let  $\delta >0$ such that \eqref{eq2martinA} holds. By setting $\bar{p}=1+\delta$ and
$\bar{q}=\frac{1+\delta}{\delta}$, we derive from the Holder
inequality that
$$F^n_A(x,\gamma)\le \frac{1}{\Gamma_n}\ES\{|\xi_1(\gamma,x,U_1)+\xi_2(\gamma,x,Z)|^{2(1+\delta)}\}^{\frac{1}{1+\delta}}
\Big(\PE(|\xi_1(\gamma,x,U_1)+\xi_2(\gamma,x,Z)|\ge
\rho\sqrt{\Gamma_n})\Big)^{\frac{\delta}{1+\delta}}$$ On the one
hand, we deduce from \eqref{elemineq} and from \eqref{eq2martinA} that
$$\ES\{|\xi_1(\gamma,x,U_1)+\xi_2(\gamma,x,Z)|^{2(1+\delta)}\}^{\frac{1}{1+\delta}}\le
C(x,\delta)(\gamma^{1+\delta}+\gamma
)^{\frac{1}{1+\delta}}\le C_1(x,\delta)\gamma^{\frac{1}{1+\delta}} $$ where $x\mapsto C_1(x,\delta)$ is locally
bounded. On the other hand, we deduce from the Chebyschev
inequality that,
\begin{align*}
\Big(\PE(|\xi_1(\gamma,x,U_1)+\xi_2(\gamma,x,Z)|\ge \rho\sqrt{\Gamma_n})\Big)^{\frac{\delta}{2+\delta}}&
\le \frac{1}{(\rho^2\Gamma_n)^{\frac{\delta}{1+\delta}}}\ES\{ |\xi_1(\gamma,x,U_1)+\xi_2(\gamma,x,Z)|^{2}\}
^{\frac{\delta}{1+\delta}}\\
&\le C_2(x,\delta,\rho)(\frac{\gamma}{\Gamma_n})^{\frac{\delta}{1+\delta}}
\end{align*}
where  $x\mapsto C_2(x,\delta,\rho)$ is locally bounded. Then, for
every  $A>0$ and $\rho>0$,
$$R^{\rho,A}_{n,1}\le C_{A,\rho}\frac{1}{\Gamma_n^{1+\frac{\delta}{1+\delta}}}
\sum_{k=1}^n\gamma_k= C_{A,\rho}\frac{1}{\Gamma_n^{\frac{\delta}{1+\delta}}}\xrightarrow{n\rightarrow+\infty}0 \qquad a.s.$$
Now, we observe $R^{\rho,A}_{n,2}$. From \eqref{contri}, we have :
$\ES\{|\xi_k^n|^2/\fek\}\le CV^{\beta}(\bar{X}_{k-1})$ with
$\beta<p/2+a-1$. Therefore,
\begin{equation*}
 R_{n,2}^{A,\rho}\le\sum_{k=1}^n 1_{\{|\bar{X}_{k-1}|\ge A\}}\ES\{|\xi_k^n|^2/\fek\}\le \sup_{|x|\ge A}\frac{V^\beta(x)}{V^{\frac{p}{2}+a-1}(x)}
 \sup_{n\in\EN}\bar{\nu}_n (V^{\frac{p}{2}+a-1})=\phi(A)\bar{\nu}_n (V^{\frac{p}{2}+a-1})
\end{equation*}
where $\phi(A)\xrightarrow{A\rightarrow+\infty}0$. Since
$\sup_{n\in\EN}\bar{\nu}_n (V^{\frac{p}{2}+a-1})<+\infty$ (see
Proposition \ref{principal}), letting
$A\rightarrow+\infty$ yields
$$R^\rho_n\xrightarrow{n\rightarrow+\infty}0\qquad a.s.\quad \forall \rho>0. $$
\end{proof}
\begin{lemme}\label{resteD} Let $a\in(0,1]$, $r\ge 0$ and $p>2$  such that $\bf{(H}_{\bf{p}}\bf{)}$,
 $\mathbf{(R_{a})}$, $\mathbf{(S_{a,r})}$ hold and $p/2+a-1>2r$.
Assume that $\ES\{U_1^{\otimes3}\}=0$ and that
 $\ES\{|U_1|^{2p}\}<+\infty$. Let $f:\ER^d\mapsto\ER$ satisfying  $\mathbf{(C_f^p)}$. Then,\\
(a) If $\Gamma_n^{(2)}/\sqrt{\Gamma_n}\rightarrow0$,
 $$\frac{1}{\sqrt{\Gamma_n}}\sum_{k=1}^n \Theta_3(\gamma_k,\bar{X}_{k-1},\bar{X}_{k,2},Z^{^{(k)}})\underset{n\rightarrow+\infty}{\overset{\PE}{\longrightarrow}}0.$$
(b) If
$\Gamma_n^{(2)}/\sqrt{\Gamma_n}\rightarrow\hat{\gamma}\in
(0,+\infty]$,
$$\frac{1}{\Gamma_n^{(2)}}\sum_{k=1}^n \Theta_3(\gamma_k,\bar{X}_{k-1},\bar{X}_{k,2},Z^{^{(k)}})\underset{n\rightarrow+\infty}{\overset{\PE}{\longrightarrow}}\int\phi_3(x)\nu(dx)$$
where $\phi_3$ is defined like in Theorem \ref{vitesseAbis}.
\end{lemme}
\noindent
The proof of this lemma  is realized in subsection \ref{preuveresteD}.
\begin{Remarque} If $Z$ is a compensated compound Poisson process,
computing the rate of convergence of $\Theta_3$ consists in
evaluating what happens after its first jump. Naturally, this
argument has no sense when the Lévy measure is not finite but the
proof and the formulation of $\phi_3$ show that it keeps some
sense in average.
\end{Remarque}
\noindent
\begin{lemme}\label{resteA} Let $a\in(0,1]$, $r\ge 0$ and $p>2$  such that $\bf{(H}_{\bf{p}}\bf{)}$,
 $\mathbf{(R_{a})}$, $\mathbf{(S_{a,r})}$ hold and $p/2+a-1>2r$.
Assume that $\ES\{U_1^{\otimes3}\}=0$ and that
 $\ES\{|U_1|^{2p}\}<+\infty$. Let $f:\ER^d\mapsto\ER$ satisfying  $\mathbf{(C_f^p)}$. Then,\\
(a) If $\Gamma_n^{(2)}/\sqrt{\Gamma_n}\rightarrow0$,
\begin{align*}
& \frac{1}{\sqrt{\Gamma_n}}\sum_{k=1}^n
\Theta_1(\gamma_k,\bar{X}_{k-1})+\Theta_2(\gamma_k,\bar{X}_{k-1},U_k)\underset{n\rightarrow+\infty}{\overset{\PE}{\longrightarrow}}0,\
\\&\frac{1}{\sqrt{\Gamma_n}}\sum_{k=1}^n
R_1(\gamma_k,\bar{X}_{k-1},U_k)+R_2(\gamma_k,\bar{X}_{k-1},U_k)+R_3(\gamma_k,\bar{X}_{k-1},\bar{X}_{k,2},Z^{^{(k)}})\underset{n\rightarrow+\infty}{\overset{\PE}{\longrightarrow}}0.
\end{align*}
(b)  If
$\Gamma_n^{(2)}/\sqrt{\Gamma_n}\xrightarrow{n\rightarrow+\infty}\hat{\gamma}\in
(0,+\infty]$, we have
\begin{align*}
 \frac{1}{\Gamma_n^{(2)}}&\sum_{k=1}^n \Theta_1(\gamma_k,\bar{X}_{k-1})\underset{n\rightarrow+\infty}{\overset{\PE}{\longrightarrow}}m_1\quad\textnormal{and}\quad
\frac{1}{\Gamma_n^{(2)}}\sum_{k=1}^n\Theta_2(\gamma_k,\bar{X}_{k-1},U_k)\underset{n\rightarrow+\infty}{\overset{\PE}{\longrightarrow}}m_2,\\
\textnormal{with}\quad &m_1=\frac  1 2 \int {D^2}f(x)b(x)^{\otimes2}\nu(dx)\\
\textnormal{and,}\quad&m_2=\int \int\frac 1 6 D^3 f(x); b(x) ;(\sigma(x)u)^{\otimes 2}+\frac{1}{24} D^4 f(x)(\sigma(x)u)^{\otimes 4}\PE_{U_1}(du)\nu(dx).
\end{align*}
At last,
$$\frac{1}{\Gamma_n^{(2)}}\sum_{k=1}^n R_1(\gamma_k,\bar{X}_{k-1},U_k)+R_2(\gamma_k,\bar{X}_{k-1},U_k)+R_3(\gamma_k,\bar{X}_{k-1},\bar{X}_{k,2},Z^{^{(k)}})\underset{n\rightarrow+\infty}{\overset{\PE}{\longrightarrow}}0.$$
\end{lemme}
\begin{proof}[\textnormal{\textbf{Proof}}]  The arguments of this
proof are quite similar to those of the previous lemma. Then, we leave it
to the reader.
\end{proof}
\subsection{Synthesis and proof of Theorem \ref{vitesseAbis}}\label{preuvevitesseA}
$\bullet$ Proof of Theorem \ref{vitesseAbis} when $\Gamma_n^{(2)}/\sqrt{\Gamma_n}\rightarrow0$: Looking into the
decomposition of $\bar{\nu}_n(Af)$ introduced in lemma
\ref{lemmeZ}, we deduce from lemmas \ref{martinA}, \ref{resteD}$(a)$
and \ref{resteA}$(a)$ that
\begin{equation}\label{gbi}
\sqrt{\Gamma_n}\bar{\nu}_n(Af)-\Big(\frac{f(\bar{X}_n)-f(\bar{X}_0)}{\sqrt{\Gamma_n}}\Big)\underset{n\rightarrow+\infty}{\xrightarrow{\cal
L}}{\cal N}\Big(0,\hat{\sigma}_f^2 \Big).
\end{equation}
Now, $f\le CV$. Then, by Proposition \ref{principal}$(a)$.\textit{iii} and Jensen's inequality,
$\ES\{f(\bar{X}_n)\}\le \ES\{V^{p}(\bar{X}_n)\}\le C\Gamma_n^\frac{1}{p}$. It
 implies that
$$\frac{f(\bar{X}_n)-f(\bar{X}_0)}{\sqrt{\Gamma_n}}\underset{n\rightarrow+\infty}{\xrightarrow{\;L^1\;}}0$$
and Theorem \ref{vitesseAbis} is obvious.\\
\noindent
$\bullet$ Proof of Theorem
\ref{vitesseAbis}  when $\Gamma_n^{(2)}/\sqrt{\Gamma_n}\rightarrow\hat{\gamma}\in
(0,+\infty]$: in this case,  $\sqrt{\Gamma_n}\le
C\Gamma_n^{(2)}$. It implies that
$$\frac{f(\bar{X}_n)-f(\bar{X}_0)}{\Gamma_n^{(2)}}\underset{n\rightarrow+\infty}{\xrightarrow{\;L^1\;}}0.$$
According to Lemmas \ref{martinA}, \ref{resteD}$(b)$ and
\ref{resteA}$(b)$, we have
\begin{equation*}
\frac{\Gamma_n}{\Gamma_n^{(2)}}\bar{\nu}_n(Af)-\Big(\frac{f(\bar{X}_n)-f(\bar{X}_0)}{\Gamma_n^{(2)}}\Big)
\begin{cases}\underset{n\rightarrow+\infty}{\xrightarrow{\;\PE\;}}m&\textnormal{ if $\hat{\gamma}=+\infty$}\\
\underset{n\rightarrow+\infty}{\xrightarrow{\cal L}}{\cal
N}\big(\hat{\gamma}m,\hat{\sigma}_f^2 \big)&\textnormal{ if
$\hat{\gamma}<+\infty$}
\end{cases}
\end{equation*}
and the result follows.
\subsection{Proof of Lemma \ref{resteD}}\label{preuveresteD}
In the proof of Lemma \ref{resteD}, we usually need to show that
some sequences tend to 0 in probability. The arguments used for this are collected in the  following
 lemma (these arguments  also work for the proof
 of Lemma \ref{resteA}).
\begin{lemme} \label{techvit} Let $a\in(0,1]$, $r\ge 0$ and $p>2$. Assume $\bf{(H}_{\bf{p}}\bf{)}$,
 $\mathbf{(R_{a})}$ and $\mathbf{(S_{a,r})}$.
Suppose that
 $\ES\{|U_1|^{2p}\}<+\infty$ and let  $(F_k)$ be a sequence of random variables
  such that $F_k$ is ${\cal F}_k$-measurable.\\
\noindent
(a)  Assume that $\Gamma_n^{(2)}/\sqrt{\Gamma_n}\rightarrow0$.\\
\noindent
i. If $|F_k|\le C\gamma_k^2V^{\frac{p}{2}+a-1}(\bar{X}_{k-1})$, then,
$1/\sqrt{\Gamma_n}\sum_{k=1}^n F_{k-1}\underset{n\rightarrow+\infty}{\overset{\PE}{\longrightarrow}}0.$\\
\noindent
\textit {ii.} If $\ES\{F_k/\fek\}=0$ and
$\ES\{|F_k|^2/\fek\}\le
C\big(\gamma_k^3V^{p}(\bar{X}_{k-1})+\gamma_k^2V^{\frac{\epsilon
p}{2}}(\bar{X}_{k-1})\big)$ with $\epsilon\in[0,1)$
then,
$1/\sqrt{\Gamma_n}\sum_{k=1}^n F_{k}\underset{n\rightarrow+\infty}{\overset{\PE}{\longrightarrow}}0$.\\
\noindent
(b) Assume
that $\Gamma_n^{(2)}/\sqrt{\Gamma_n}\xrightarrow{n\rightarrow+\infty}\hat{\gamma}\in
(0,+\infty]$. Then,\\
\noindent
i. If $|F_k| \le C\gamma_k^{2+\delta}V^{\frac{p}{2}+a-1}(\bar{X}_{k-1})$,
$1/\Gamma_n^{(2)}\sum_{k=1}^n F_{k-1}\underset{n\rightarrow+\infty}{\overset{\PE}{\longrightarrow}}0.$\\
\noindent
ii. If $\ES\{F_k/\fek\}=0$ and $\ES\{|F_k|^2/\fek\}\le
C\big(\gamma_k^3V^{\epsilon p}(\bar{X}_{k-1})+\gamma_k^2V^{\frac{\epsilon
p}{2}}(\bar{X}_{k-1})\big)$ with $\epsilon\in[0,1)$ then,
$1/\Gamma_n^{(2)}\sum_{k=1}^n F_{k}\underset{n\rightarrow+\infty}{\overset{\PE}{\longrightarrow}}0.$
\end{lemme}
\begin{proof}[\textnormal{\textbf{Proof}}] \textit{(a)} i.  By Proposition \ref{principal}$(a).iii$,
$\ES\{V^p(\bar{X}_n)\}\le C\Gamma_n$. We then derive from Jensen's inequality that
\begin{equation*}
\frac{1}{\sqrt{\Gamma_n}}\ES\{\sum_{k=1}^n |F_{k-1}|\}\le
\frac{1}{\sqrt{\Gamma_n}}\sum_{k=1}^n\gamma_k^2\Gamma_k^\frac{\bar{p}}{p}.
\end{equation*}
where $\bar{p}=p/2+a-1$. Hence, the first assertion is obvious if
\begin{equation}\label{oteh}
\frac{1}{\sqrt{\Gamma_n}}\sum_{k=1}^n\gamma_k^2\Gamma_k^\frac{\bar{p}}{p}\xrightarrow{n\rightarrow+\infty}0.
\end{equation}
\noindent
If \eqref{oteh} is not fulfilled, then we have $\liminf \frac{1}{\sqrt{\Gamma_n}}\sum_{k=1}^n\gamma_k^2\sqrt{\Gamma_k}>0$
because $\bar{p}/p\le 1/2$. It follows from the Kronecker Lemma that we have necessary $\sum_{k\ge 1}
\gamma_k^2=+\infty$.
By setting $\eta_k=\gamma_k^2$, we can apply
Proposition \ref{principal} and deduce that
\begin{equation}\label{suptendu}
\sup_{n\ge 1}\frac{1}{\Gamma_n^{(2)}}\sum_{k=1}^n \gamma_k^2 V^{\frac{p}{2}+a-1}(\bar{X}_{k-1})<+\infty\quad a.s.
\end{equation}
Since
$\Gamma_n^{(2)}/\sqrt{\Gamma_n}\xrightarrow{n\rightarrow+\infty}0$,
the first assertion follows when \eqref{oteh} is not fulfilled.\\
\noindent
ii. Since  $\ES\{V^p(\bar{X}_n)\}\le \Gamma_n$, we derive from Jensen's inequality that  $\ES\{|F_k|^2\}\le
C(\gamma_k^3 \Gamma_k+\gamma_k^2 \Gamma_k^{\frac{\epsilon }{2}}),$
$\textnormal{with $\epsilon \in [0,1)$.}$ On the one hand, one checks that
$\sum_{k=1}^n \gamma_k^3\Gamma_k\le (\Gamma_n^{(2)})^2$. Hence, since
$\Gamma_n^{(2)}/\sqrt{\Gamma_n}\xrightarrow{n\rightarrow+\infty}0$,
we have
\begin{equation}\label{ret2}
\frac{1}{\Gamma_n}\sum_{k=1}^n\gamma_k^3 \Gamma_k\le
C\frac{(\Gamma_n^{(2)})^2}{(\sqrt{\Gamma_n})^2}\xrightarrow{n\rightarrow+\infty}0.
\end{equation}
On the other hand, one observes that $\sum_{k=1}^n\frac{\gamma_k^2
\Gamma_k^{\frac{\epsilon }{2}}}{\Gamma_k}\le
\sum_{k=1}^n\frac{\gamma_k^2}{(\Gamma_k^{(2)})^{2-\epsilon}}<+\infty.$
Hence, the  Kronecker Lemma implies that
\begin{equation}\label{ret}
\frac{1}{\Gamma_n}\sum_{k=1}^n\gamma_k^2 \Gamma_k^{\frac{\epsilon
}{2}}\xrightarrow{n\rightarrow+\infty}0.
\end{equation}
It follows that
$\frac{1}{\Gamma_n}\sum_{k=1}^n\ES\{|F_k|^2\}\xrightarrow{n\rightarrow+\infty}0$.
This yields the second assertion of $(a)$.\\
\noindent \textit{(b)} $i.$  We derive from the assumptions that
\begin{equation}\label{1a}
\frac{1}{\Gamma_n^{(2)}}\sum_{k=1}^n |F_k|\le
\frac{C}{\Gamma_n^{(2)}}
 \sum_{k=1}^n \gamma_k^{2+\delta} V^{\frac{p}{2}+a-1}(\bar{X}_{k-1}).
 \end{equation}
Then, $(b).i.$ follows from \eqref{suptendu} which
is still valid because $\sum_{k\ge1}\gamma_k^2=+\infty$.\\
\noindent
$ii.$  It
suffices to check that
\begin{equation}\label{1c}
\frac{1}{(\Gamma_n^{(2)})^2}\sum_{k=1}^n\ES\{|F_k|^2\}\xrightarrow{n\rightarrow+\infty}+\infty.
\end{equation}
With the same arguments as in $(a).ii$, one checks that
\begin{equation}\label{1489}
\frac{1}{(\Gamma_n^{(2)})^2}\sum_{k=1}^n
\ES\{|F_k|^2\}\le\frac{C}{(\Gamma_n^{(2)})^2}\sum_{k=1}^n\gamma_k^2(\Gamma_k^{(2)})^{\epsilon}+
\frac{C}{(\Gamma_n^{(2)})^2}\sum_{k=1}^n\gamma_k^3\Gamma_k^{\epsilon}\quad\textnormal{with
$\epsilon\in[0,1)$.} \end{equation} On the one hand,  we deduce
from the Kronecker Lemma that
$$\frac{1}{(\Gamma_n^{(2)})^2}\sum_{k=1}^n\gamma_k^2(\Gamma_k^{(2)})^{\epsilon}\xrn{n\nrn}0.$$
On the other hand, for every $\epsilon\in[0,1)$
$$\frac{1}{(\Gamma_n^{(2)})^2}\sum_{k=1}^n\gamma_k^3\Gamma_k^{\epsilon}\overset{n\rightarrow+\infty}{=}
\frac{1}{(\Gamma_n^{(2)})^2}\sum_{k=1}^n\gamma_k^3\Gamma_k^{\epsilon}<+\infty$$
because $\sum_{k=1}^n \gamma_k^3\Gamma_k\le (\Gamma_n^{(2)})^2$. Hence, we derive  from \eqref{1489} that
$$\frac{1}{(\Gamma_n^{(2)})^2}\sum_{k=1}^n
\ES\{|F_k|^2\}\xrn{n\nrn}0$$ which completes the proof.
\end{proof}

\noindent
\textbf{Proof of Lemma \ref{resteD}.} $(a)$ In order to alleviate the notations, we prove
the lemma in the one-dimensional case. We set
$$\bar{\Theta}_3(\gamma,x,z,Z)=\int_0^\gamma ds\int\pi(dy_1)
\Big(\tilde{H}^f(z+\kappa(x) Z_s,x,y_1)-
\tilde{H}^f(x,x,y_1)\Big).$$ $\bar{\Theta}_3(\gamma,x,z,Z)$ is the
compensator of $\Theta_3(\gamma,x,z,Z)$. Then, since
$\Theta_3(\gamma,x,z,Z)$ is a purely  discontinuous process, we
have
$$\ES\{|(\Theta_3-\bar{\Theta}_3)(\gamma,x,z,Z)|^2\}
=\ES\{\int_0^\gamma ds\int\pi(dy_1)\big|\tilde{H}^f(z+\kappa(x) Z_s,x,y_1)- \tilde{H}^f(x,x,y_1)\big|^2\}.$$
By Taylor's formula,
$$\tilde{H}^f(z,x,y_1)=\int_0^1 \left( f'(z+\theta\kappa(x) y_1)-f'(x)\right)\kappa(x)y_1d\theta.$$
It follows that,
\begin{align*}
|\tilde{H}^f(z+\kappa(x) Z_s,x,y_1)-& \tilde{H}^f(x,x,y_1)|
\\&\le
\sup_{\theta\in [0,1]}
|f'(z+\kappa(x)( Z_s+\theta y_1))-f'(x+\theta\kappa(x) y_1)|.|\kappa(x)y_1|\\&
+|f'(z+ \kappa(x)Z_s)-f'(x)|.|\kappa(x)y_1|.
\end{align*}
$f'$ is a Lipschitz continuous function. Then, by setting $z=x+\gamma
b(x)+\sqrt{\gamma}\sigma(x)u$, we deduce from the assumptions on
the coefficients and from the fact that $\ES\{|Z_s|^2\}=O(s)$ that
\begin{align}
\ES\{|(\Theta_3-&\bar{\Theta}_3)(\gamma,x,z,Z)|^2\}\nonumber\\
&\le
C\ES\{\int_0^\gamma ds\int\pi(dy_1)\Big(\gamma^2 |b|^2(x)+\gamma
|\sigma(x)|^2u^2+|\kappa(x)|^2|Z_s|^2\Big)|\kappa(x)|^2
y_1^2\}\nonumber\\&\le C\big(\gamma^3V^{a+r}(x)+\gamma^2(1+|u|^2)V^{2r}(x)\big).\label{390}
\end{align}
Set
$F_k=(\Theta_3-\bar{\Theta}_3)(\gamma_k,\bar{X}_{k-1},\bar{X}_{k,2},Z^{^{(k)}})$.
Since $\ES\{F_k/\fek\}=0$, $a+r<p$ and $2r<p/2+a-1$,  it follows from Lemma
\ref{techvit}$(a).ii$ and from the preceding inequality that
\begin{equation}\label{etap1}
\frac{1}{\sqrt{\Gamma_n}}\sum_{k=1}^n
(\Theta_3-\bar{\Theta}_3)(\gamma_k,\bar{X}_{k-1},\bar{X}_{k,2},Z^{^{(k)}})\xrightarrow{\PE}0\quad \textnormal{when }n\rightarrow+\infty.
\end{equation}
Now, since $f^{(2)}$ is bounded, we deduce from Taylor's formula that
\begin{equation*}
|\tilde{H}^f(z+\kappa(x) Z_s,x,y_1)- \tilde{H}^f(x,x,y_1)|\le C\|\kappa(x)\|^2.|y_1|^2.
\end{equation*}
Then,
$$\ES\{|\bar{\Theta}_{3}(\gamma,x,z,Z)|^2\}\le C\gamma^2\|\kappa(x)\|^4\int|y_1|^2\pi(dy_1)^2\le C\gamma^2 V^{2r}(x).$$
By Lemma \ref{techvit}$(a).ii$, it follows that,
\begin{equation}\label{etap2}
\frac{1}{\sqrt{\Gamma_n}}\sum_{k=1}^n
\Big(\bar{\Theta}_{3}(\gamma_k,\bar{X}_{k-1},\bar{X}_{k,2},Z^{^{(k)}})-\ES\{\bar{\Theta}_{3}(\gamma_k,\bar{X}_{k-1},\bar{X}_{k,2},Z^{^{(k)}})/\fek\}\Big)\underset{\PE}{\xrightarrow{n\rightarrow+\infty}}0.
\end{equation}
Then, by \eqref{etap1} and \eqref{etap2}, $(a)$ is
obvious if we prove that
\begin{align}
&\frac{1}{\sqrt{\Gamma_n}}\sum_{k=1}^n
\ES\{\bar{\Theta}_{3,1}(\gamma_k,\bar{X}_{k-1},\bar{X}_{k,2},Z^{^{(k)}})/\fek\}\underset{n\rightarrow+\infty}{\overset{\PE}{\longrightarrow}}0\qquad\textnormal{and,}\label{puix}\\&
\frac{1}{\sqrt{\Gamma_n}}\sum_{k=1}^n
\ES\{\bar{\Theta}_{3,2}(\gamma_k,\bar{X}_{k-1},\bar{X}_{k,2},Z^{^{(k)}})/\fek\}\underset{n\rightarrow+\infty}{\overset{\PE}{\longrightarrow}}0\qquad
\textnormal{with}\label{puiz}
\end{align}
\begin{align*}
&\bar{\Theta}_{3,1}(\gamma,x,z,Z)\\ &=\gamma \bar{\pi}\int_0^1
ds\int\tilde{\pi}(dy_1)\int_0^1d\theta \left(
f^{(2)}(z+\kappa(x)(Z_{s\gamma}+\theta y_1))
-f^{(2)}(z+\theta \kappa(x)y_1)\right)(1-\theta)\kappa^2(x),\\
&\bar{\Theta}_{3,2}(\gamma,x,z)=\gamma \int\tilde{\pi}(dy_1)
\int_0^1d\theta\left( f^{(2)}(z+\theta\kappa(x)
y_1)-f^{(2)}(x+\theta \kappa(x)y_1)\right) (1-\theta)\kappa^2(x)
\end{align*}
where $\bar{\pi}= \int y_1^2\pi(dy) (<+\infty)$ and $\tilde{\pi}$ is a probability measure
defined by $\tilde{\pi}(dy_1)=y_1^2\pi(dy_1)/\bar{\pi}.$
Let us prove \eqref{puix}. By Itô's formula, we have
\begin{align*}
 f^{(2)}(z+\kappa(x)(Z_{s\gamma}+\theta y_1))-f^{(2)}(z+\theta \kappa(x)y_1)
 &=\int_0^{s\gamma}f^{(3)}\big(z+\kappa(x)(Z_{v^{-}}+\theta y_1)\big)\kappa(x) dZ_v\\&
+\sum_{0<v<s\gamma}\tilde{H}^{f^{(2)}}\big(z+\kappa(x)(Z_{v^-}+\theta
y_1),x,\Delta Z_v\big).\end{align*} Since $f^{(3)}$ is bounded, the
first term of the right-hand side is a martingale. Therefore,  we obtain by the
compensation formula that
\begin{align*}
 &\ES\{\bar{\Theta}_{3,1}(\gamma,x,z,Z)\}\\&=\gamma \bar{\pi}\int_0^1
ds\int\tilde{\pi}(dy_1)\int_0^1d\theta\int_0^{s\gamma}
dv\int\pi(dy_2)
\ES\big\{\tilde{H}^{f^{(2)}}(z+\kappa(x)(Z_{v}+\theta y_1),x,
y_2)\big\}(1-\theta)\kappa^2(x).
\end{align*}
Finally, since
$$|\tilde{H}^{f^{(2)}}\big(z+\kappa(x)(Z_{v}+\theta y_1),x, y_2\big)|\le C\|f^{(4)}\|_\infty\kappa^2(x)y_2^2,$$
we have
$$\ES\{|\bar{\Theta}_{3,1}(\gamma_k,\bar{X}_{k-1},\bar{X}_{k,2},Z^{^{(k)}})|/ \fek\}
\le C\gamma_k^2 \kappa^4(\bar{X}_{k-1})\le C\gamma_k^2 V^{2r}(\bar{X}_{k-1}).$$ Since $2r\le p/2+a-1$,
 \eqref{puix} follows from  Lemma
\ref{techvit}$(a).i$.\\
\noindent
Now, let us prove  \eqref{puiz}. Set
$z=x+\gamma b(x)+\sqrt{\gamma}\sigma(x)U_1$. By Taylor's
formula, there exist $$\xi_1 \in [
x+\sqrt{\gamma}\sigma(x)U_1+\theta \kappa(x)y_1, z+\theta
\kappa(x)y_1] \quad \textnormal{and}\quad \xi_2 \in [x+\theta
\kappa(x)y_1, x+\sqrt{\gamma}\sigma(x)U_1+\theta \kappa(x)y_1]$$
such that
\begin{align*}
f^{(2)}(z+\theta&\kappa(x) y_1)-f^{(2)}(x+\kappa(x)\theta
\kappa(x)y_1)\\&=\gamma f^{(3)}(\xi_2)b(x)
+\sqrt{\gamma}f^{(3)}(x+\theta \kappa(x)y_1)\sigma(x)U_1+\gamma
f^{(4)}(\xi_1)\sigma^2(x)U_1^2.
\end{align*}
Since   $f^{(3)}$ and $f^{(4)}$ are bounded and $U_1$ is centered,
$$|\ES\{\bar{\Theta}_{3,2}(\gamma_k,\bar{X}_{k-1},\bar{X}_{k,2})/\fek\}|\le
 C\gamma_k^2\big(|b|\kappa^2(\bar{X}_{k-1})+\sigma^2\kappa^2(\bar{X}_{k-1})\big)\le
  C\gamma_k^2  V^{{a\vee 2r}}(\bar{X}_{k-1})$$
where we have used that $a/2+r\le a\vee 2r$. Since $a\vee 2r\le p/2+a-1$,
we deduce \eqref{puiz} from   Lemma \ref{techvit}$(a)$.i.\\
\\
\noindent $b)$ We keep the  notations of $(a)$. On the one hand, by \eqref{390}
and Lemma \ref{techvit}$(b).ii$, we have
\begin{equation*}
\frac{1}{\Gamma_n^{(2)}}\sum_{k=1}^n
(\Theta_3-\bar{\Theta}_3)(\gamma_k,\bar{X}_{k-1},\bar{X}_{k,2},Z^{^{(k)}})
\underset{\PE}{\xrightarrow{n\rightarrow+\infty}}0.
\end{equation*}
On the other hand, we have
$$\ES\{|\bar{\Theta}_3(\gamma,x,z,Z)|^2\}\le C\gamma^2 V^{\frac{\epsilon p}{2}}(x)$$
with $\epsilon\in[0,1)$. Hence, by
 Lemma \ref{techvit}$(b).ii$, it follows that
\begin{equation*}
\frac{1}{\Gamma_n^{(2)}}\sum_{k=1}^n
\Big(\bar{\Theta}_{3}(\gamma_k,\bar{X}_{k-1},\bar{X}_{k,2},Z^{^{(k)}})-\ES\{\bar{\Theta}_{3}(\gamma_k,\bar{X}_{k-1},\bar{X}_{k,2},Z^{^{(k)}})/\fek\}\Big)\underset{\PE}{\xrightarrow{n\rightarrow+\infty}}0.
\end{equation*}
Finally, it suffices to prove that
\begin{align}
&\frac{1}{\Gamma_n^{(2)}}\sum_{k=1}^n
\ES\{\bar{\Theta}_{3,1}(\gamma_k,\bar{X}_{k-1},\bar{X}_{k,2},Z^{^{(k)}})/\fek\}
\underset{n\rightarrow+\infty}{\overset{\PE}{\longrightarrow}}m_{3,1}\label{puix2}\\
\textnormal{and,}\quad&\frac{1}{\Gamma_n^{(2)}}\sum_{k=1}^n
\ES\{\bar{\Theta}_{3,2}(\gamma_k,\bar{X}_{k-1},\bar{X}_{k,2},Z^{^{(k)}})/\fek\}\underset{n\rightarrow+\infty}{\overset{\PE}{\longrightarrow}}m_{3,2}.\label{puiz2}
\end{align}
with $m_{3,1}+m_{3,2}=\int \phi_3(x)\nu(dx)$. In order to prove
\eqref{puix2}, we first show that
\begin{equation}\label{toure}
\frac{1}{\Gamma_n^{(2)}}\sum_{k=1}^n\Big(\ES\{\bar{\Theta}_{3,1}(\gamma_k,\bar{X}_{k-1},\bar{X}_{k,2},Z^{^{(k)}})
/\fek\}-\bar{\Theta}_{3,1}(\gamma_k,\bar{X}_{k-1},\bar{X}_{k-1},0)\Big)\underset{n\rightarrow+\infty}{\overset{\PE}{\longrightarrow}}0.
\end{equation}
Since $f^{(4)}$ is a bounded and Lipschitz continuous  function, $f^{(4)}$ is also $2\delta$-Holder
for every $\delta\in (0,1/2]$, $i.e.$
$$[f^{(4)}]_{2\delta}=\sup_{x,y\in\ER^d}\frac{|f^{(4)}(y)-f^{(4)}(x)|}{|y-x|^{2\delta}}<+\infty.$$
It follows from the Taylor formula that
\begin{align*}
\Big|\tilde{H}^{f^{(2)}}(z+\kappa(x)(Z_{v}+\theta y_1),x,
y_2)-&\tilde{H}^{f^{(2)}}(z+\kappa(x)(Z_{v}+\theta y_1),x,
y_2)\Big|\\&\le
C[f^{(4)}]_{2\delta}\Big(|z-x|^{2\delta}+\kappa(x)^{2\delta}|Z_v|^{2\delta}\Big)\kappa(x)^2y_2^2.
\end{align*}
By setting $z=x+\gamma b(x)+\sqrt{\gamma}\sigma(x)u$ and taking
$\delta$ sufficiently small, we have
$$\ES\{|\tilde{H}^{f^{(2)}}(z+\kappa(x)(Z_{v}+\theta y_1),x, y_2)
-\tilde{H}^{f^{(2)}}(x+\kappa(x)\theta y_1,x, y_2)|\}\le
C\gamma^\delta (1+|u|^{2\delta})V^{\frac{p}{2}+a-1}(x).$$
This implies that
\begin{align*}
&\big|\ES\Big\{\bar{\Theta}_{3,1}(\gamma_k,\bar{X}_{k-1},\bar{X}_{k,2},Z^{^{(k)}})
-\bar{\Theta}_{3,1}(\gamma_k,\bar{X}_{k-1},\bar{X}_{k-1},0)/\fek\Big\}\big|\le
C\gamma_k^{2+\delta}V^{\frac{p}{2}+a-1}(\bar{X}_{k-1}).
\end{align*}
Then, \eqref{toure} follows from Lemma \ref{techvit}$(b).i$. \\
\noindent Now,  $\bar{\Theta}_{3,1}(\gamma_k,\bar{X}_{k-1},\bar{X}_{k-1},0)$
is $\fek$-measurable and
$$|\bar{\Theta}_{3,1}(\gamma_k,\bar{X}_{k-1},\bar{X}_{k-1},0)|\le C\gamma_k^2
V^{2r}(\bar{X}_{k-1}).$$ Since $2r<p/2+a-1$,  we can apply Proposition
\ref{principal}  with $\eta_k=\gamma_k^2$. We obtain
\begin{align*}
\frac{1}{\Gamma_n^{(2)}}\sum_{k=1}^n \bar{\Theta}_{3,1}(\gamma_k,&\bar{X}_{k-1},\bar{X}_{k-1},0)
=\frac{1}{\Gamma_n^{(2)}}\sum_{k=1}^n\frac{\gamma_k^2}{2}\phi_{3,1}(\bar{X}_{k-1})\xrightarrow{n\rightarrow+\infty} \int \phi_{3,1}(x)\nu(dx)\quad a.s.\\
\textnormal{with}\quad
\phi_{3,1}(x)&=\frac{1}{2}\int\pi(dy_1)\int_0^1d\theta\int\pi(dy_2)
\tilde{H}^{f^{(2)}}\big(x+\kappa(x)\theta y_1,x,
y_2\big)(1-\theta)\kappa^2(x)y_1^2\\&=\frac{1}{2}\int\pi(dy_1)\int\pi(dy_2)\tilde{H}^{\tilde{H}^f_{.,x,y_1}}(x,x,y_2).
\end{align*}
It follows from \eqref{toure} that
$$\frac{1}{\Gamma_n^{(2)}}\sum_{k=1}^n\bar{\Theta}_{3,1}
(\gamma_k,\bar{X}_{k-1},\bar{X}_{k,2},Z^{^{(k)}})/\fek\}\xrightarrow{n\rightarrow+\infty}
 m_{3,1}=\int\phi_{3,1}(x)\nu(dx)\quad a.s.$$
Finally, we prove \eqref{puiz2}. Since $f^{(3)}$ and $f^{(4)}$ are
bounded and Lipschitz continuous, we deduce from   Taylor's formula that for
every $\delta \in[0,1/2]$,
\begin{align*}
&\ES\{\bar{\Theta}_{3,2}(\gamma_k,\bar{X}_{k-1},\bar{X}_{k,2})/\fek\}=\gamma_k^2(\phi_{3,2}(\bar{X}_{k-1})+\phi_{3,3}(\bar{X}_{k-1}))+ \bar{\rho}_1(\bar{X}_{k-1},\gamma_k)+\bar{\rho}_2(\bar{X}_{k-1},\gamma_k)\\
 & \textnormal{with}\quad\phi_{3,2}(x)=\int\pi(dy_1)\int_0^1d\theta f^{(3)}(x+\theta \kappa(x)y_1)b(x)(1-\theta)(\kappa(x)y_1)^2\\
&\phi_{3,3}(x)=\int\pi(dy_1)\int_0^1d\theta\int\PE_{U_1}(du)
f^{(4)}(x+\theta \kappa(x)y_1)\sigma^2(x)
u^2(1-\theta)(\kappa(x)y_1)^2\\& |\bar{\rho}_1(x,\gamma)|\le
[f^{(3)}]_\delta \gamma^{2+\delta}
|b(x)|^{1+\delta}|\kappa(x)|^2\quad\textnormal{and}\quad
|\bar{\rho}_2(x,\gamma)|\le [f^{(4)}]_{2\delta} \gamma^{2+\delta}
|\sigma(x)|^{2(1+\delta)}|\kappa(x)|^2.
\end{align*}
Since $(a/2+r)\vee (2r)<p/2+a-1$, one can find  $\delta>0$ such that
 $$|b(x)|^{1+\delta}|\kappa(x)|^2+|\sigma(x)|^{2(1+\delta)}|\kappa(x)|^2\le V^{\frac{p}{2}+a-1}(x).$$
On the one hand,  Lemma  \ref{techvit}$(b).i$ and
the assumptions on
 the coefficients allow us to conclude that
$$\frac{1}{\Gamma_n^{(2)}}\sum_{k=1}^n\Big(\bar{\rho}_1(\bar{X}_{k-1},\gamma_k)+\bar{\rho}_2(\bar{X}_{k-1},\gamma_k)\Big)\xrightarrow{n\rightarrow+\infty}0.$$
On the other hand, as
$|b|(x)|\kappa|^2(x)+|\sigma|^2(x)|\kappa|^2(x)=o(V^{\frac{p}{2}+a-1}(x))$ when $|x|\rightarrow+\infty$, we derive from Proposition
\ref{principal} applied with $\eta_k=\gamma_k^2$  that
$$ \frac{1}{\Gamma_n^{(2)}}\sum_{k=1}^n\gamma_k^2(\phi_{3,2}(\bar{X}_{k-1})+\phi_{3,3}(\bar{X}_{k-1}))\xrightarrow{n\rightarrow+\infty}m_{3,2}$$
with $m_{3,2}=\int (\phi_{3,2}(x)+\phi_{3,3}(x))\nu(dx)$. Checking
that
$$\phi_{3,2}(x)+\phi_{3,3}(x)=\int\pi(dy_1)\Big((H^f_{.,x,y_1})'(x)b(x)+\int\PE_{U_1}(du)(H^f_{.,x,y_1})^{(2)}(x)(\sigma(x)u)^2\Big).$$
completes the proof.
\section{\large{Proof of Theorem \ref{vitesseWbis}}}\label{proofvitesseW}
The  proof is built as follows. Like in the proof of Theorem \ref{vitesseAbis}, we firstly decompose $\bar{\nu}^{^P}_n(Af)$ and
$\bar{\nu}^{^W}_n(Af)$ (see Lemma
\ref{lemmeZB}).  Some new terms appear due to the approximation of the jump component. That is why in the sequel, we
focus on these parts of the decomposition (see Lemmas
\ref{resteB} and \ref{resteBtilde2}). The other terms can be
studied by the same process as their corresponding terms in the
decomposition of $\bar{\nu}^{}_n(Af)$ and then, are left to the
reader. \noindent We denote by $(Z^{^{(k)_P}})_{k\ge1}$, a
sequence of independent and càdlàg processes such that
$$(Z_t^{^{(k)_P}})_{t\ge 0}\overset{\cal L}{=}(Z_{t,k})_{t\ge 0}\quad\textnormal{and}
\quad Z_{\gamma_k}^{^{(k)_P}}=\bar{Z}_k^{^P}\quad \forall k\ge 1$$
For a ${\cal C}^2$-function $f$ such that ${D^2}f$ is bounded, we  define $A^{k,P}$ and $A^{k,W}$ by
\begin{align*}
& A^{k,P}f(x)=\psg\nabla f,b\psd(x)+\frac{1}{2}{\rm{Tr}}(\sigma^* {D^2}
f\sigma)(x)+\int_{\{|y|>u_k\}} \tilde{H}^f(x,x,y)\pi(dy)\\&
A^{k,W}f(x)=A^{k,P}f(x)+\frac{1}{2}\int_{\{|y|\le u_k\}}
{D^2}f(x)(\kappa(x)y)^{\otimes2}\pi(dy),
\end{align*}
These operators correspond respectively to the infinitesimal generators of
\begin{align*}
&dX_{t}=b(X_{t^{-}})dt+\sigma(X_{t^{-}})dW_{t}+\kappa(X_{t^{-}})dZ_{t,k}\quad\textnormal{and}
\quad\\&
dX_{t}=b(X_{t^{-}})dt+\sigma(X_{t^{-}})dW_{t}+\kappa(X_{t^{-}})d(Z_{t,k}+Q_k
\tilde{W}_t),
\end{align*}
 where $\tilde{W}$ is a $q$-dimensional Brownian motion independent
of $W$ and $(Z_{t,k})_{t\ge 0}$.
\begin{lemme} \label{lemmeZB} For a $C^2$-function $f$ such that ${D^2}f$ is bounded, we have
the following decompositions
\begin{align*}
1)\quad\sum_{k=1}^n\gamma_k Af(\bar{X}^{^P}_{k-1})&=G_n^{^P}
+f(\bar{X}^{^P}_n)-f(x)-\sum_{k=1}^n \Big(\xi_1(\gamma_k,\bar{X}^{^P}_{k-1},U_k)+\xi^k_{2,P}(\gamma_k,\bar{X}^{^P}_{k-1},Z^{^{(k)_P}})\Big)\\
&-\sum_{k=1}^n
\Big(\Theta_1(\gamma_k,\bar{X}^{^P}_{k-1})+\Theta_2(\gamma_k,\bar{X}^{^P}_{k-1},U_k)+\Theta_3(\gamma_k,\bar{X}^{^P}_{k-1},\bar{X}_{k,2}^{^P},Z^{^{(k)_P}})\Big)\\&
-\sum_{k=1}^n
\Big((R_1+R_2)(\gamma_k,\bar{X}^{^P}_{k-1},U_k)+R_3(\gamma_k,\bar{X}^{^P}_{k-1},\bar{X}_{k,2}^{^P},Z^{^{(k)_P}})\Big),
\end{align*}
where $\bar{X}_{k,2}^{^P}=\bar{X}^{^{P}}_{k-1}+\gamma_k
b(\bar{X}^{^{P}}_{k-1})+\sqrt{\gamma_k}\sigma(\bar{X}^{^{P}}_{k-1})U_k,$
$G_n^{^P}=\sum_{k=1}^n \gamma_k (Af-A^{k,P}
f)(\bar{X}^{^P}_{k-1})$ and for a càdlàg process $Y$,
$$\xi_{2,P}^k(\gamma,x,Y)=\int_0^\gamma\psg\nabla f(x),\kappa(x)dY_s\psd+\Big(\sum_{0<s\le \gamma}\tilde{H}^f(x,x,\Delta Y_s))-\gamma\int_{D_k}\tilde{H}^f(x,x,y)\pi(dy)\Big).$$
\begin{align*}
2) \sum_{k=1}^n\gamma_k Af(\bar{X}^{^{W}}_{k-1})&=G_n^{^W}+J_n^{^W}+f(\bar{X}^{^{W}}_n)-f(x)-
\sum_{k=1}^n \Big(\xi_1(\gamma_k,\bar{X}^{^{W}}_{k-1},U_k)+\xi^k_{2,B}
(\gamma_k,\bar{X}^{^{W}}_{k-1},Z^{^{(k)_P}})\Big)\\
&-\sum_{k=1}^n
\Big(\Theta_1(\gamma_k,\bar{X}^{^{W}}_{k-1})+\Theta_2(\gamma_k,\bar{X}^{^{W}}_{k-1},U_k)
+\Theta_3(\gamma_k,\bar{X}^{^{W}}_{k-1},\bar{X}_{k,2}^{^P},Z^{^{(k)_P}})\Big)\\&
-\sum_{k=1}^n
\Big((R_1+R_2)(\gamma_k,\bar{X}^{^{W}}_{k-1},U_k)+R_3(\gamma_k,\bar{X}^{^{W}}_{k-1},
\bar{X}_{k,2}^{^P},Z^{^{(k)_P}})\Big),
\end{align*}
where
$\bar{X}^{^{W}}_{k,2}=\bar{X}^{^{W}}_{k-1}+\bar{X}^{^{W}}_{k-1}+\gamma_k
b(\bar{X}^{^{W}}_{k-1})+\sqrt{\gamma_k}\sigma(\bar{X}^{^{W}}_{k-1})U_k,$
$$G_n^{^W}=\sum_{k=1}^n \gamma_k
(Af-A^{k,W}f)(\bar{X}^{^{W}}_{k-1}),\; J_n^{^W}=-\sum_{k=1}^n
f(\bar{X}^{^{W}}_{k})-f(\bar{X}^{^{W}}_{k,3})+\gamma_k({A}^{k,P} f-A^{k,W}
f)(\bar{X}^{^{W}}_{k-1}),$$ with
$\bar{X}^{^{W}}_{k,3}=\bar{X}^{^{W}}_k-\kappa(\bar{X}^{^{W}}_{k-1})Q_k\Lambda_k.$
\end{lemme}
\noindent
The two following lemmas are devoted to the additional terms of
the preceding decomposition. In Lemma \ref{resteB}, we compute the
rate of $G_n^{^P}$ and $G_n^{^W}$ and in Lemma \ref{resteBtilde2},
we show that $J_n^{^W}$ does not have any consequences on the rate
of the procedure.
\begin{lemme} \label{resteB} Let $a\in(0,1]$, $r\ge 0$ and $p>2$  such that $\bf{(H}_{\bf{p}}\bf{)}$,
 $\mathbf{(R_{a})}$, $\mathbf{(S_{a,r})}$ hold and $2r<p/2+a-1$.
Suppose  that
 $\ES\{|U_1|^{2p}\}<+\infty$ and $\ES\{|\Lambda_1|^{2p}\}<+\infty$.
 Let $f:\ER^d\mapsto\ER$ satisfying  $\mathbf{(C_f^p)}$. Then,\\
 \noindent
(1) i. If $\lim_{n\rightarrow+\infty}\beta_{n,\pi}^{(2)}<+\infty$,
$\frac{1}{\sqrt{\Gamma_n}}\sum_{k=1}^n \gamma_k (Af-A^{k,P}
f)(\bar{X}^{^P}_{k-1})\underset{n\rightarrow+\infty}{\overset{\PE}{\longrightarrow}}0$.\\
\noindent ii.  If
$\lim_{n\rightarrow+\infty}\beta_{n,\pi}^{(2)}=+\infty$,
$$\limsup_{n\rightarrow+\infty}\frac{1}{\beta_{n,\pi}^{(2)}}\sum_{k=1}^n \gamma_k
\big|(Af-A^{k,P} f)(\bar{X}^{^P}_{k-1})\big|\le \bar{m}_2
 \quad a.s.$$
\noindent
where $\bar{m}_2=\frac{\|{D^2}f\|_{\infty}}{2}\int\|\kappa\|^2(x)\nu(dx)$. Furthermore, if $\mathbf{(A_2^2)}$ holds,
\begin{equation}\label{convP1}
\frac{1}{\beta_{n,\pi}^{(2)}}\sum_{k=1}^n \gamma_k
(Af-A^{k,P} f)(\bar{X}^{^P}_{k-1})\xrightarrow{n\rightarrow+\infty}m_2\quad a.s.\quad\textnormal{with }
|m_2|\le \bar{m}_2.
\end{equation}
 (2) Assume that $s=3$ or that $s=4$ if $\pi$ is  quasi-symmetric in the neighborhood of 0.\\
 i. If $\lim_{n\rightarrow+\infty}\beta_{n,\pi}^{(s)}<+\infty$
$\frac{1}{\sqrt{\Gamma_n}}\sum_{k=1}^n \gamma_k (Af-A^{k,W} f)(\bar{X}^{^{W}}_{k-1})\underset{n\rightarrow+\infty}{\overset{\PE}{\longrightarrow}}0$.\\
ii. If $\lim_{n\rightarrow+\infty}\beta_{n,\pi}^{(s)}=+\infty$,
$$\limsup_{n\rightarrow+\infty}\frac{1}{\beta_{n,\pi}^{(3)}}\sum_{k=1}^n \gamma_k
\Big|(Af-A^{k,W} f)(\bar{X}^{^{W}}_{k-1})\Big|\le \bar{m}_{s}
\quad a.s. $$ where
$\bar{m}_{s}=C_s\|D^sf\|_{\infty}\int\|\kappa(x)\|^s\nu(dx)$ with
$C_3=\frac{d^{\frac 3 2}}{6}$ and $ C_4=\frac{{d^2}}{24}$. Furthermore, if
, $\mathbf{(A_s^2)}$ holds,
\begin{equation}\label{convW1}
\frac{1}{\beta_{n,\pi}^{(s)}}\sum_{k=1}^n \gamma_k
(Af-A^{k,W} f)(\bar{X}^{^W}_{k-1})\xrightarrow{n\rightarrow+\infty}m_s\quad a.s.\quad\textnormal{ with }
|m_s|\le \bar{m}_s.
\end{equation}
\end{lemme}
\begin{proof}[\textnormal{\textbf{Proof}}]
$(1)i.$ By  Taylor's formula, we have
\begin{equation*}
|Af(x)-A^{k,P}
f(x)|=\frac{1}{2}\int {D^2}f(\xi_y)(\kappa(x)y)^{\otimes2}\pi(dy)\le\frac{\|{D^2}f\|_{\infty}}{2}
\int\sum_{i,j}|(\kappa(x)y)_i(\kappa(x)y)_j|\pi(dy).
\end{equation*}
For $z\in\ER^d$, $|\sum_{i,j}z_iz_j|\le |z|_1^2\le d|z|^2$. It follows that
\begin{equation}\label{1k}
|Af(\bar{X}^{^P}_{k-1})-A^{k,P}
f(\bar{X}^{^P}_{k-1})|\le\frac{d}{2}\|{D^2}f\|_{\infty}
\|\kappa(\bar{X}^{^P}_{k-1})\|^2\int_{|y|\le
u_k}|y|^2\pi(dy).
\end{equation}
Since $r\le p/2$ and $\ES\{V^p(
\bar{X}^{^P}_{k-1})\}\le \Gamma_k$, we deduce that
\begin{equation*} \sum_{k=1}^n\gamma_k\ES\{|Af(\bar{X}^{^P}_{k-1})-A^{k,P}
f(\bar{X}^{^P}_{k-1})|\}\le \sum_{k=1}^n \gamma_k\int_{|y|\le
u_k}|y|^2\pi(dy)\sqrt{\Gamma_k}.
\end{equation*}
Now, as $\lim_{n\rightarrow+\infty}\beta_{n,\pi}^{(2)}<+\infty$,
 the Kronecker Lemma yields
$$\frac{1}{\sqrt{\Gamma_n}}\sum_{k=1}^n \gamma_k\int_{|y|\le u_k}|y|^2\pi(dy)\sqrt{\Gamma_k}\xrightarrow{n\rightarrow+\infty}0.$$
The first assertion is obvious.\\
\noindent $ii.$ Since
$\lim_{n\rightarrow+\infty}\beta_{n,\pi}^{(2)}=+\infty$, we deduce
from Proposition \ref{principal}  with
$\eta_k=\gamma_k\int_{|y|\le u_k}|y|^2\pi(dy)$ that
$$\frac{1}{\beta_{n,\pi}^{(2)}}\sum_{k=1}^n\gamma_k\int_{|y|\le u_k}|y|^2\pi(dy)
\|\kappa(\bar{X}^{^P}_{k-1})\|^2\xrightarrow{n\rightarrow+\infty}\int\|\kappa(x)\|^2\nu(dx)
\quad a.s.$$ because $\|\kappa\|^2=o(V^{\frac{p}{2}+a-1})$. Then,
the  second assertion follows from \eqref{1k}. \\
\noindent Assume now that $\mathbf{(A_2^2)}$ holds. Since ${D^2}f$ is Lipschitz continuous,
 we deduce from Taylor's formula that
\begin{align*}
&Af(x)-A^{k,P}
f(x)=\frac{1}{2}\Big(\sum_{i,j}\rho_k(i,j)\psi_{i,j}(x)+R^k_{i,j}(x)\Big)\\
&\textnormal{with}\quad \rho_k(i,j)=\int_{\{|y|\le u_k\}}y_iy_j\pi(dy),
\quad\psi_{i,j}(x)=\sum_{l,m}\kappa_{i,l}\frac{\partial^2 f}{\partial_l\partial_m}\kappa_{m,j}(x))
\end{align*}
and $|R^k_{i,j}(x)|\le \int_{\{|y|\le u_k\}}|y|^3\pi(dy)\|\kappa(x)\|^3.$\\
\noindent
According to $\mathbf{(A_2^2)}$, for every $i,j$,
$\lim\rho_k(i,j)/\int_{\{|y|\le u_k\}}|y|^2\pi(dy)=\alpha_{i,j}\in\ER.$
Set $\eta_k=\gamma_k\int_{\{|y|\le u_k\}}|y|^2\pi(dy)$ and $H_n=\sum_{k=1}^n\eta_k$. Then,
\begin{align*}
\frac{1}{\beta_{n,\pi}^{(2)}}\sum_{k=1}^n \gamma_k
\rho_k(i,j)\psi_{i,j}(\bar{X}^{^P}_{k-1})=
\frac{\alpha_{i,j}}{H_n}
\sum_{k=1}^n\eta_k \psi_{i,j}(\bar{X}^{^P}_{k-1})+\frac{1}{H_n}
\sum_{k=1}^n\varepsilon^1_k\eta_k\psi_{i,j}(\bar{X}^{^P}_{k-1})
\end{align*}
with $\varepsilon^1_k=(\rho_k(i,j)-\alpha_{i,j}\eta_k)/\eta_k$. Firstly, since $\psi_{i,j}\le CV^r$ and $r<p/2+a-1$,
 Proposition \ref{principal} applied with $\eta_k=\gamma_k\int_{\{|y|\le u_k\}}|y|^2\pi(dy)$ yields
$$\frac{1}{H_n}\sum_{k=1}^n
\eta_k \psi_{i,j}(\bar{X}^{^P}_{k-1})\xrightarrow{n\rightarrow+\infty}
\int \psi_{i,j}(x)\nu(dx)\quad a.s.$$
Secondly, $\varepsilon^1_k=o(1)$. Since
$\sup_{n\ge1}1/H_n\sum_{k=1}^n\eta_k V^{\frac{p}{2}+a-1}(\bar{X}^{^P}_{k-1})<+\infty$ $a.s.$, it is then easy to check that
$$\frac{1}{H_n}
\sum_{k=1}^n\varepsilon^1_k\eta_k\psi_{i,j}(\bar{X}^{^P}_{k-1})\xrightarrow{n\rightarrow+\infty}0\quad a.s.$$
The same argument is appropriate for  $R^k_{i,j}$ with
 $\varepsilon^2_k=(\int_{\{|y|\le u_k\}}|y|^3\pi(dy))/(\int_{\{|y|\le u_k\}}|y|^2\pi(dy))$.
Finally, we obtain
\begin{equation*}
\frac{1}{\beta_{n,\pi}^{(2)}}\sum_{k=1}^n \gamma_k
(Af-A^{k,P} f)(\bar{X}^{^P}_{k-1})\xrightarrow{n\rightarrow+\infty}m_2=\sum_{i,j}\frac{\alpha_{i,j}}{2}\int\psi_{i,j}(x)\nu(dx).
\end{equation*}

\noindent 2) We derive from  Taylor's formula
\begin{equation}\label{1m}
|Af(x)-A^{k,W}f(x)|\le
C_s\|D^sf\|_{\infty}\|\kappa(x)\|^s\int_{|y|\le u_k}|y|^s\pi(dy).
\end{equation}
with $s=3$ and $C_3=\frac{d^{\frac 3 2}}{6}$, or $s=4$ and
$C_4=\frac{d^{2}}{24}$  if $\int_{|y|\le u_k}
y^{\otimes3}\pi(dy)=0$. As $2r\le p/2$ and $\ES\{V^p(
\bar{X}^{^{W}}_{k-1})\}\le \Gamma_k$, it implies that
$$\sum_{k=1}^n\gamma_k\ES\{|Af-A^{k,W}
f|(\bar{X}^{^{W}}_{k-1})\}\le  C\sum_{k=1}^n \gamma_k\int_{|y|\le
u_k}|y|^s\pi(dy)\sqrt{\Gamma_k},$$ with $s=3$ or $s=4$ if
$\int_{|y|\le u_k} y^{\otimes3}\pi(dy)=0$. If
$\lim_{n\rightarrow+\infty}\beta_{n,\pi}^{(s)}<+\infty$, we derive
from the  Kronecker Lemma that
$$\frac{1}{\sqrt{\Gamma_n}}\sum_{k=1}^n \gamma_k\int_{|y|\le u_k}|y|^s\pi(dy)\sqrt{\Gamma_k}\xrightarrow{n\rightarrow+\infty}0.$$
and  the first assertion of (2) follows.\\
Assume now that $\lim_{n\rightarrow+\infty}\beta_{n,\pi}^{(s)}=+\infty$.
Applying Proposition \ref{principal} to $f(x)=\|\kappa(x)\|^s$
with $\eta_k=\gamma_k\int_{|y|\le u_k}|y|^s\pi(dy)$ yields
$$\limsup_{n\rightarrow+\infty}\frac{1}{\beta_{n,\pi}^{(s)}}\sum_{k=1}^n\gamma_k\int_{|y|\le
u_k}|y|^s\pi(dy)\|\kappa(\bar{X}^{^{W}}_{k-1})\|^s<+\infty\quad a.s.$$
Then,  $(2).ii$ follows from \eqref{1m}.\\
\noindent Finally, the proof of \eqref{convW1} is similar to that of \eqref{convP1}.
\end{proof}
\begin{lemme}\label{resteBtilde2}Let $a\in(0,1]$, $r\ge 0$ and $p>2$  such that $\bf{(H}_{\bf{p}}\bf{)}$,
 $\mathbf{(R_{a})}$, $\mathbf{(S_{a,r})}$ hold and $2r<p/2+a-1$.
Suppose  that
 $\ES\{|U_1|^{2p}\}<+\infty$ and that $\ES\{|\Lambda_1|^{2p}\}<+\infty$.
 Let $f:\ER^d\mapsto\ER$ satisfying  $\mathbf{(C_f^p)}$. Then,\\
\begin{equation}
\frac{1}{\sqrt{\Gamma_n\vee \Gamma_n^{(2)}}}\sum_{k=1}^n
\Big(f(\bar{X}^{^{W}}_{k})-f(\bar{X}^{^{W}}_{k,3})+\gamma_k({A}^{k,P} f-A^{k,W}
f)(\bar{X}^{^{W}}_{k-1})\Big)\underset{n\rightarrow+\infty}{\overset{\PE}{\longrightarrow}}0.
\end{equation}
\end{lemme}
\begin{proof}[\textnormal{\textbf{Proof}}]
Since
$$({A}^{k,W} f-A^{k,P}f)(x)=\frac{1}{2}\int_{\{|y|\le u_k\}}
{D^2}f(x)(\kappa(x)y)^{\otimes2}\pi(dy)=\frac{1}{2}\ES\{{D^2} f(x) (\kappa(x)Q_k\Lambda_k)^{\otimes2}\},$$
we derive from  Taylor's formula that
\begin{align*}
f(\bar{X}^{^{W}}_{k})-f(\bar{X}^{^{W}}_{k,3})&+\gamma_k({A}^{k,P} f-A^{k,W}) f(\bar{X}^{^{W}}_{k-1})=
\xi_1(\gamma_k,\bar{X}^{^{W}}_{k,3},\bar{X}^{^{W}}_{k-1},Q_k\Lambda_k)\\&+\tilde{R}_{k,1}(\gamma_k,\bar{X}^{^{W}}_{k-1}, Q_k\Lambda_k)
+ \tilde{R}_{k,2}(\gamma_k,\bar{X}^{^{W}}_{k,3},\bar{X}^{^{W}}_{k-1}, Q_k\Lambda_k)
\end{align*}
\begin{align*}
&\textnormal{where,}\quad\xi_1(\gamma,z,x,Q_k\Lambda_k)=\sqrt{\gamma}\psg\nabla f(z),\kappa(x)Q_k\Lambda_k\psd ,\\
&\tilde{R}_{k,1}(\gamma,x,Q_k\Lambda_k)=\frac {\gamma}{2} \big( {D^2}
f(x) (\kappa(x)Q_k\Lambda_k)^{\otimes2}-
\ES\{{D^2} f(x) (\kappa(x)Q_k\Lambda_k)^{\otimes2}\}\big)\\
&\tilde{R}_{k,2}(\gamma,z,x,Q_k\Lambda_k)=\gamma\int_0^1
\big({D^2}f(z+\theta\sqrt{\gamma}\kappa(x)Q_k\Lambda_k)-{D^2}f(x)\big)(1-\theta)
(\kappa(x)Q_k\Lambda_k)^{\otimes2}d\theta.
\end{align*}
Setting $\theta_n=\sqrt{\Gamma_n}\vee \Gamma_n^{(2)}$, it suffices
to show the three following steps:\\

\noindent a)
$\theta_n^{-1}\sum_{k=1}^n\xi_{1}(\gamma_k,\bar{X}^{^{W}}_{k-1},Q_k\Lambda_k)
\underset{n\rightarrow+\infty}{\overset{\PE}{\longrightarrow}}0$,\\
b) $\theta_n^{-1}\sum_{k=1}^n
\tilde{R}_{k,1}(\gamma_k,\bar{X}^{^{W}}_{k-1},Q_k\Lambda_k)\underset{n\rightarrow+\infty}{\overset{\PE}{\longrightarrow}}0$.\\
c) $\theta_n^{-1}\sum_{k=1}^n
\tilde{R}_{k,2}(\gamma_k,\bar{X}^{^{W}}_{k,3},\bar{X}^{^{W}}_{k-1},Q_k\Lambda_k)\underset{n\rightarrow+\infty}{\overset{\PE}{\longrightarrow}}0.$
\\

\noindent a)  We set \begin{align*}
&\xi_1(\gamma,z,x,Q_k\Lambda_k)=\xi_{1,1}(\gamma,x,Q_k\Lambda_k)+\xi_{1,2}(\gamma,z,x,Q_k\Lambda_k)
\end{align*}
with $\xi_{1,1}(\gamma,x,v)=\sqrt{\gamma}\psg\nabla
f(x),\kappa(x)v\psd $ and
$\xi_{1,2}(\gamma,z,x,v)=\sqrt{\gamma}\psg\nabla
f(z)-\nabla f(x),\kappa(x)v\psd $.
Let $(M_{n,1})$ and
$(M_{n,2})$ be the $({\cal F}_n)$-martingales defined  by
$$M_{n,1}=\sum_{k=1}^n\xi_{1,1}(\gamma_k,\bar{X}^{^{W}}_{k-1},Q_k\Lambda_k)\quad
\textnormal{and}\quad
M_{n,2}=\sum_{k=1}^n\xi_{1,2}(\gamma_k,\bar{X}^{^{W}}_{k,3},\bar{X}^{^{W}}_{k-1},Q_k\Lambda_k).$$
We have to prove that $\theta_n^{-1}M_{n,1}\xrn{n\nrn}0$ and $\theta_n^{-1}M_{n,2}\xrn{n\nrn}0$.\\
\noindent According to
$\mathbf{(C_f^p)}$ and the assumptions on $\kappa$, we check that
\begin{align*}
\frac{<M>_{n,1}}{\Gamma_n}\le\frac{C}{\Gamma_n}\sum_{k=1}^n \gamma_k\int_{|y|\le
u_k}|y|^2\pi(dy)V^{\frac{p}{2}+a-1}(\bar{X}^{^{W}}_{k-1}).
\end{align*}
Now, by \eqref{Vsuptendu}, $\sup_{n\ge 1}1/\Gamma_n\sum_{k=1}^n
\gamma_k V^{\frac{p}{2}+a-1}(\bar{X}^{^{W}}_{k-1})<+\infty$ $a.s.$ Since $\int_{|y|\le
u_k}|y|^2\pi(dy)\rightarrow0$,
that $<M>_{n,1}/\Gamma_n\xrn{n\nrn}0$ $a.s$. Then,
\begin{equation}\label{xiun}
\frac{1}{\sqrt{\Gamma_n}\vee \Gamma_n^{(2)}}|M_{n,1}|\le\frac{1}{\sqrt{\Gamma_n}}|M_{n,1}|\underset{n\rightarrow+\infty}{\overset{\PE}{\longrightarrow}}0.
\end{equation}
Now, we turn to $(M_{n,2})$.
Since $\nabla f$ is Lipschitz continuous, it follows from the assumptions on the coefficients that
\begin{align*}
&\ES\{|\xi_{1,2}(\gamma_k,\bar{X}^{^{W}}_{k,3},\bar{X}^{^{W}}_{k-1},Q_k\Lambda_k)|^2/\fek\}\le
C\int_{|y|\le u_k}|y|^2\pi(dy)\Big(\gamma_k^3 V^{\epsilon
p}(\bar{X}^{^{W}}_{k-1})+\gamma_k^2
V^{\epsilon\frac{p}{2}}(\bar{X}^{^{W}}_{k-1})\Big).
\end{align*}
with $\epsilon<1$ and $\bar{p}=p/2+a-1$.
Then, by a variant of Lemma \ref{techvit}$(a).ii$ and $(b).ii$, one checks that
\begin{equation*}
\frac{<M_{n,2}>}{(\sqrt{\Gamma_n}\vee \Gamma_n^{(2)})^2}\xrightarrow{n\rightarrow+\infty}0
\quad a.s.\Longrightarrow\quad\frac{1}{\sqrt{\Gamma_n}\vee
\Gamma_n^{(2)}}M_{n,2}\underset{L^2}{\xrightarrow{n\rightarrow+\infty}}0.
\end{equation*}
\noindent b) $\tilde{R}_{k,1}$ is very closed to
$R_2$ introduced in Lemma \ref{lemmeZ} and the
arguments are similar.\\
c) As ${D^2}f$ is bounded, one observes that
$$\ES\{ |\tilde{R}_{k,2}(\gamma_k,\bar{X}^{^{W}}_{k,2},\bar{X}^{^{W}}_{k-1},Q_k\Lambda_k)|^2/\fek\}
\le C\gamma_k^2 V^{2 r}(\bar{X}^{^{W}}_{k-1}).$$ Then, since $2r<p/2$,
a variant of Lemma \ref{techvit}$(a).ii$ and $(b).ii$ yields
\begin{equation}\label{rie}
\frac{1}{\sqrt{\Gamma_n}\vee \Gamma_n^{(2)}}\sum_{k=1}^n
\Big(\tilde{R}_{k,2}(\gamma_k,\bar{X}^{^{W}}_{k,2},\bar{X}^{^{W}}_{k-1},Q_k\Lambda_k)-\ES\{
\tilde{R}_{k,2}(\gamma_k,\bar{X}^{^{W}}_{k,2},\bar{X}^{^{W}}_{k-1},Q_k\Lambda_k)/\fek\}\Big)\underset{n\rightarrow+\infty}{\overset{\PE}{\longrightarrow}}0.
\end{equation}
By setting
$\zeta_{k,\theta}(x)=x+\gamma_k b(x)+\sqrt{\gamma_k} \sigma(x)U_k+\kappa(x)(\bar{Z}_k+\theta Q_k\Lambda_k)$, we
decompose the integrand of $\tilde{R}_{k,2}$ as follows:
$${D^2}f(\zeta_{k,\theta}(x))-{D^2}f(x)=\big({D^2}f(x+\gamma_k b(x))-{D^2}f(x)\big)+
\big({D^2}f(\zeta_{k,\theta}(x))-{D^2}f(x+\gamma_k b(x))\big).$$
On the one hand, set $y_k=\sqrt{\gamma_k}Q_k\Lambda_k$. By Taylor's formula, we have
\begin{align*}
&({D^2}f(x+\gamma_k b(x))-{D^2}f(x))y_k^{\otimes2}= D^3 f(\xi_k^1);\gamma_k b(x);y_k^{\otimes2}
\end{align*}
where  $D^3 f(u);v;y^{\otimes2}:=\sum_{i,j}\psg\nabla {D^2}f_{i,j}(u),v\psd y_iy_j$ and
$\xi_k^1\in[x,x+\gamma_k b(x)]$. Thus, since $D^3f$ is bounded, we deduce that
\begin{align*}
\Big|\ES\Big\{D^3 f(\bar{X}_{k-1}^{^W});
\gamma_k b(\bar{X}_{k-1}^{^W});
(\kappa(\bar{X}_{k-1}^{^W})Q_k\Lambda_k)^{\otimes2}/{\cal F}_{k-1}^{^W}\Big\}\Big|\le C \gamma_k\int_{|y|\le
u_k}|y|^2\pi(dy)V^{\frac{a}{2}+r}(\bar{X}^{^{W}}_{k-1}).
\end{align*}
On the other hand, set $\Delta_\theta(\gamma_k,x)=\zeta_{k,\theta}(x)-(x+\gamma b(x))$. By Taylor's formula,
\begin{align*}
({D^2}f(\zeta_{k,\theta}(x))-{D^2}f(x+\gamma b(x)))y_k^{\otimes2}&=
D^3 f(x);\Delta_\theta(\gamma_k,x);y_k^{\otimes2}+\frac{1}{2} D^4f(\xi_k^2);(\Delta_\theta(\gamma_k,x))^{\otimes 2};y_k^{\otimes2}\end{align*}
where  $D^4f(u);v^{\otimes 2};y^{\otimes2}=\sum_{i,j}{D^2}({D^2}f_{i,j}(u))v^{\otimes2}y_iy_j$
and $\xi_k^2\in[x+\gamma b(x),\zeta_{k,\theta}(x)]$.
The random variables $U_k$, $\bar{Z}_k$ and $\Lambda_k$ are independent and independent of ${\cal F}_{k-1}^{^W}$. Then, since
$\ES\{U_k/{\cal F}_{k-1}^{^W}\}=
\ES\{\bar{Z}_k/{\cal F}_{k-1}^{^W}\}=\ES\{\Lambda_k^{\otimes3}/{\cal F}_{k-1}^{^W}\}=0$, we have
\begin{align*}
\ES\Big\{D^3 f(\bar{X}_{k-1}^{^W});
\Delta_\theta(\gamma_k,\bar{X}_{k-1}^{^W});
(\kappa(\bar{X}_{k-1}^{^W})Q_k\Lambda_k)^{\otimes2}/{\cal F}_{k-1}^{^W}\Big\}=0.
\end{align*}
Now, since $D^4f$ is bounded, one checks that
\begin{align*}
\ES\Big\{|D^4 f(\bar{X}_{k-1}^{^W});
\big(\Delta_\theta(\gamma_k,\bar{X}_{k-1}^{^W})\big)^{\otimes2};
(\kappa(\bar{X}_{k-1}^{^W})Q_k\Lambda_k)^{\otimes2}/{\cal F}_{k-1}^{^W}|\Big\}\le C \int_{|y|\le
u_k}|y|^2\pi(dy)\gamma_k V^{2r}(\bar{X}^{^{W}}_{k-1}).
\end{align*}
Since $a/2+r\le p/2+a-1$ and $2r\le p/2+a-1$, it follows that
\begin{align*}
|\ES\{
\tilde{R}_{k,2}(\gamma_k,\bar{X}^{^{W}}_{k,2},\bar{X}^{^{W}}_{k-1},Q_k\Lambda_k)/\fek\}|&\le
C\int_{|y|\le
u_k}|y|^2\pi(dy)
\gamma_k^2V^{\frac{p}{2}+a-1}(\bar{X}^{^{W}}_{k-1})
\end{align*}
By a variant of Lemma \ref{techvit}$(a).i$ and $(b).i$, we derive from the previous inequality
 that
$$\frac{1}{\sqrt{\Gamma_n}\vee\Gamma_n^{(2)}}\sum_{k=1}^n
\ES\{ \tilde{R}_{k,2}(\gamma_k,\bar{X}^{^{W}}_{k,2},\bar{X}^{^{W}}_{k-1},Q_k\Lambda_k)/\fek\}
\underset{n\rightarrow+\infty}{\xrightarrow{\PE}}0,$$
Then, assertion $c)$ follows from \eqref{rie}.
\end{proof}
\section{\large{An additional result}}\label{extending}
In this section, we present a partial extension when
 the Lévy process has a moment of order $2p$ with $p\in(1,2]$.
In this case, stating some global results as in Theorems \ref{vitesseAbis} and \ref{vitesseWbis} would need
two kinds
of restrictions: either to assume that at least the derivatives of $f$ tend to 0 when $|x|\rightarrow+\infty$
or to impose  more constraints on the growth of the coefficients. The first alternative
leads to a very technical proof and the second one can not be really envisaged for the drift term. Actually,
we recall that in this type of problem, $b$ produces the mean-reverting effect and then, it would  not be  natural
to suppose that for instance, $b$ is bounded .\\
That is why we propose  to state   a partial result for fast-decreasing steps for which
the extension does only require some weak restrictions on $f$.  We introduce
a new assumption on the steps depending
on the intensity of the mean-reverting:
\begin{equation}\label{condpasextension}
 \frac{\Gamma_n^{(2)}}{\sqrt{\Gamma_n}}
\xrightarrow{n\rightarrow+\infty}0\quad\textnormal{if $a=1$}\quad\textnormal{and,}\quad
\frac{1}{\sqrt{\Gamma_n}}\sum_{k=1}^n
\gamma_k^2\Gamma_k^{\frac{a\vee(2r)}{p}}\xrightarrow{n\rightarrow+\infty}0\quad\textnormal{if
$a<1.$}
\end{equation}
Then,
\begin{theorem}\label{vitesseA}
Assume that $\ES\{|Z_t|^{2p}\}<+\infty$ with $p\in(1,2]$ and that \eqref{edss} admits a
unique
 invariant measure $\nu$ .
 Let $a\in(0,1]$ and $r\ge 0$
such that  $\mathbf{(R_{a})}$  and
$\mathbf{(S_{a,r})}$ are satisfied and such that $p/2+a-1>r$. Suppose that
$\ES\{U_1^{\otimes3}\}=0$, $\ES\{|U_1|^{4}\}<+\infty$ and $\eta_n=\gamma_n$ for every $n\ge 1$.
Let $f:\ER^d\mapsto\ER$ be a ${\cal C}^4$-function  having bounded derivatives and satisfying
$f(x)=O(\sqrt{V(x)})$ as $|x|\rightarrow+\infty$. Then,\\
(a) Scheme (E): If \eqref{condpasextension} holds, $\sqrt{\Gamma_n}\bar{\nu}_n(Af)\underset{n\rightarrow+\infty}{\xrightarrow{\cal
L}}{\cal N}\Big(0,\hat{\sigma}_f^2 \Big).$\\
(b) Scheme (P): If \eqref{condpasextension} holds and $\beta_{n,\pi}^{(2)}/\sqrt{\Gamma_n}\rightarrow0,$
the conclusion of $(a)$ is valid for Scheme (P).\\
(c) Scheme (W): If \eqref{condpasextension} holds and $\beta_{n,\pi}^{(s)}/\sqrt{\Gamma_n}\rightarrow0,$ with
$s=3$ or $s=4$ if $\pi$ is quasi-symmetric in the neighborhood
of 0,
the conclusion of $(a)$ is valid for Scheme (W).
\end{theorem}
\begin{Remarque}  We refer to \cite{bib26} for a proof of this result. \\
\noindent Assumption \eqref{condpasextension} is less constraining when $a=1$ because
the $L^p$-control of the Euler scheme is better in this case (see Proposition \ref{principal}).
Note that when $a=1$, Theorem \ref{vitesseA}(a) corresponds to Theorem \ref{vitesseAbis}(a) when
$\hat{\gamma}=0$. \\
\noindent
Let $\gamma_k=\gamma_1k^{-\zeta}$ with $\zeta\in (0,1]$ and $\gamma_1>0$. For Scheme (E),
Theorem \ref{vitesseA} applies in the following cases:
$$\zeta>\frac{1}{3}\quad\textnormal{if $a=1$}\quad\textnormal{and}\quad \zeta>\frac{p+2\eta}{3p+2\eta}
\quad \textnormal{if $a<1$}$$ where $\eta=a\vee(2r)$. Since
$\sqrt{\Gamma_n}\overset{\rightarrow+\infty}{\sim}\sqrt{\frac{\gamma_1}{1-\zeta}}n^{\frac{1-\zeta}{2}}$
if $\zeta\in(0,1)$, we derive   that for every $\epsilon>0$, there
exists an Euler scheme with polynomial step such that the rate of
convergence is of order $n^{\frac{1}{3}-\epsilon}$ if $a=1$ and
$n^{\frac{p}{3p+2\eta}-\epsilon}$ if $a<1$.\\
\noindent
\end{Remarque}
\section{\large{Numerical comparison of Schemes (P) and (W)}}\label{simulations}
\begin{figure}[h]
\begin{minipage}[c]{.5\linewidth}\centering
\includegraphics[width=5cm]{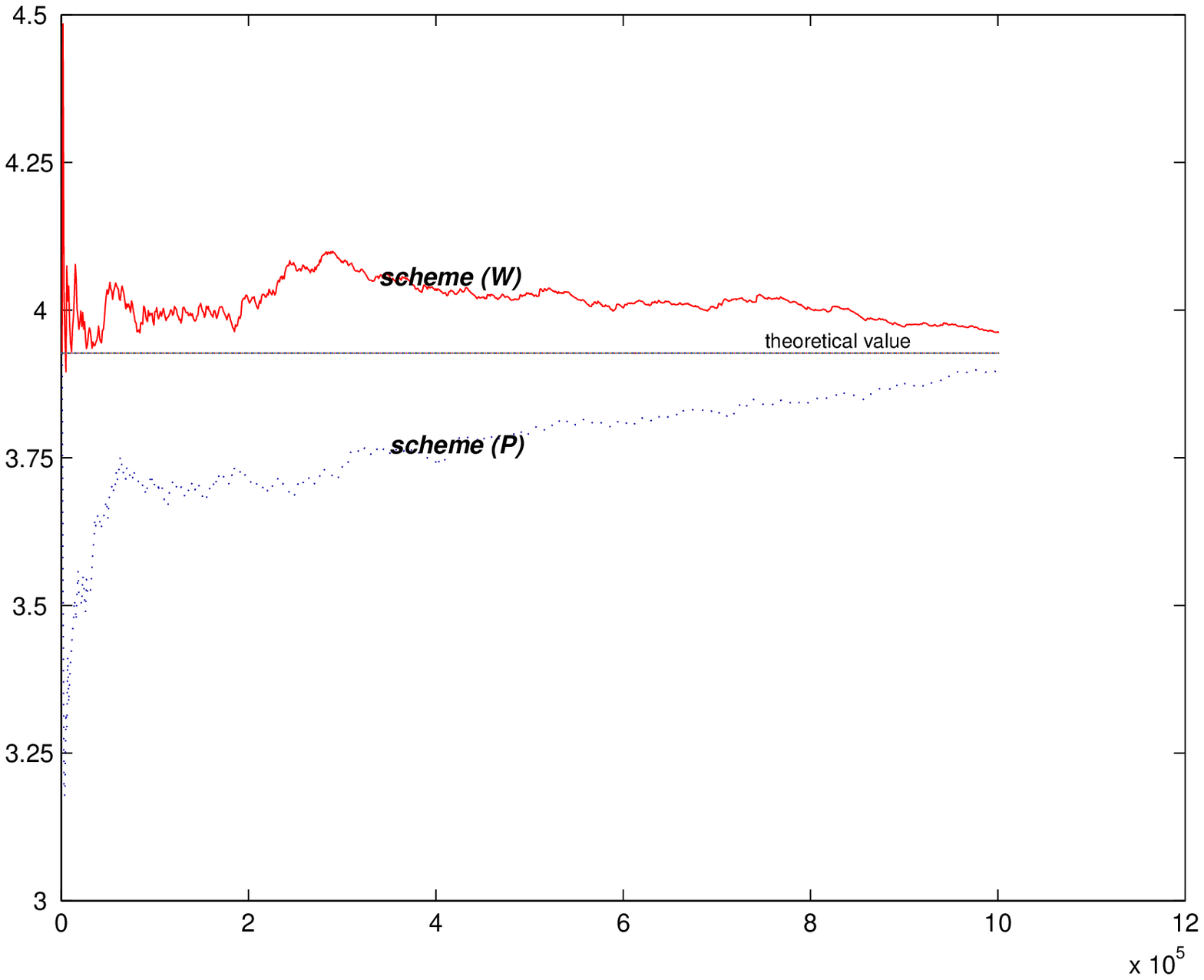}
\caption{$n\mapsto \bar{\nu}_n(\phi)$, $\alpha=1$}\label{figsim1}
\end{minipage}\hfill
\begin{minipage}[c]{.5\linewidth}\centering\label{figsim2}
\includegraphics[width=5cm]{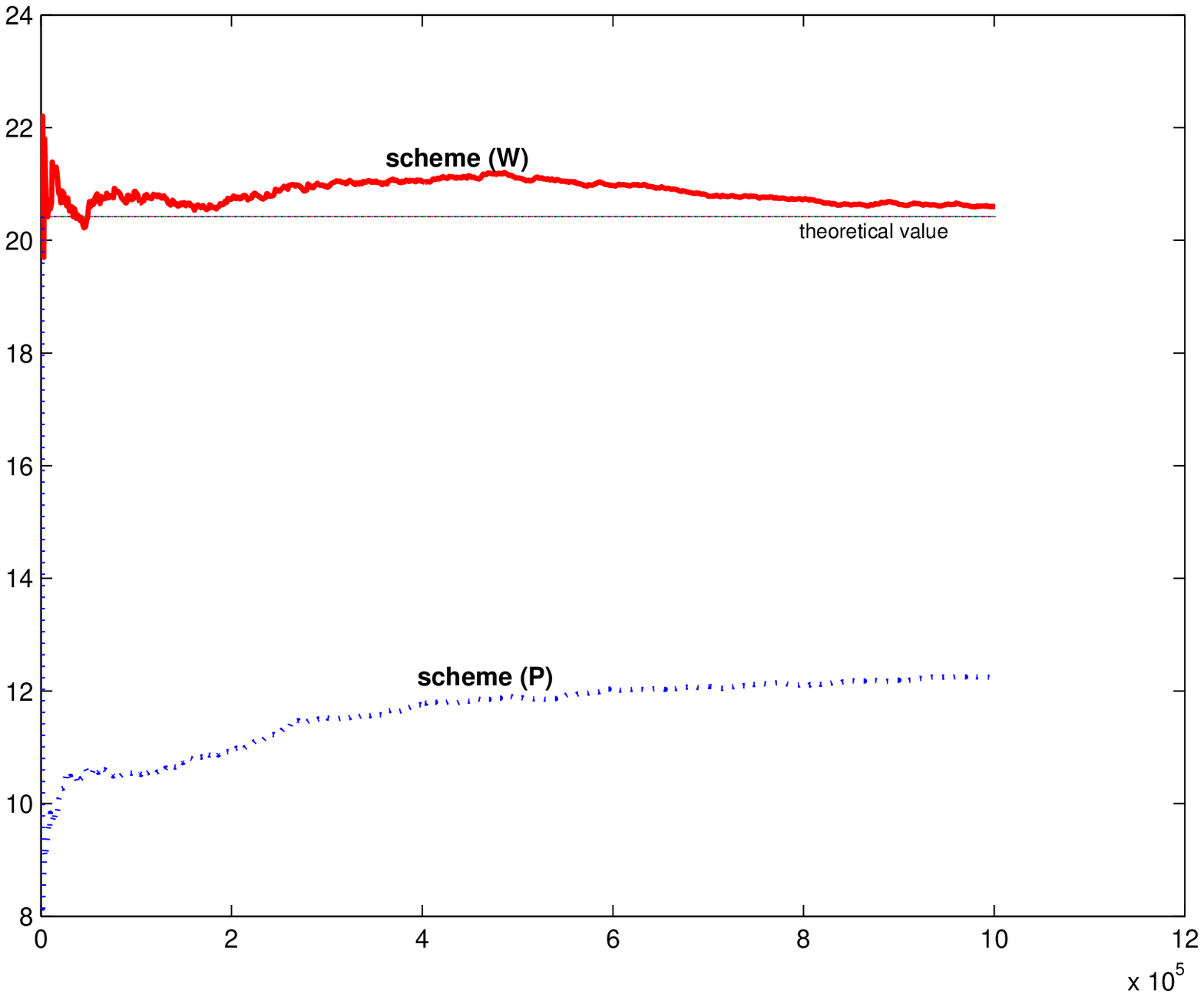}
\caption{$n\mapsto \bar{\nu}_n(\phi)$, $\alpha=5/3$}
\end{minipage}
\end{figure}

\textbf{1. When $\nu(f)$ can be theoretically computed.} In this
first example, we are interested in the two-dimensional  SDE
\begin{equation}\label{OUR}
dX_t=-X_{t^-}dt+dZ_t
\end{equation}
where $Z$ is a    symmetric purely discontinuous Lévy process
(having no drift term). We consider $\phi:x\mapsto |x|^2$ and
denote by $\nu$, the unique  invariant SDE \eqref{OUR}. We can
easily compute  $\nu(\phi)$. In fact, as  $\pi$ is symmetric
 and $\nu(A\phi)=0$,
$$ A\phi=-2\phi+\int |y|^2\pi(dy)\quad\Longrightarrow\quad \nu(\phi)=\frac{1}{2}\int |y|^2\pi(dy).$$
Let us test this theoretical result on a $2$-dimensional example.
Assume that
$$\pi^{(\alpha)}(dy)=1_{\{|y|\le 1\}}\frac{1}{|y|^{\alpha+2}}\lambda_2(dy)
+1_{\{|y|> 1\}}\frac{1}{|y|^8}\lambda_2(dy)\quad \textnormal{with }\alpha\in(1,3).$$
We have $\nu(\phi)=\pi(1/(2-\alpha)+1/4)$. In  figures
\ref{figsim1} and \ref{figsim2}, we observe the rate for two
values of $\alpha$ taking the choices of steps and truncation
thresholds of Proposition \ref{vitpart}$(b)$.
 These simulations are coherent with the theoretical results. Indeed, when
$\alpha=1$, the optimal asymptotic  rates induced by Schemes (P)
and (W) are the same (with order $n^{\frac{1}{3}}$). When
$\alpha=5/3$,  the optimal asymptotic rate induced by scheme
(W) is still of order $n^{\frac{1}{3}}$ whereas that of  scheme
(P) is of order $n^{\frac{1}{7}}$.\\
\noindent \textbf{2. Another example.} Now, we observe the
following two-dimensional SDE
\begin{equation*}
dX_t=-\frac{X_{t^-}}{\sqrt{1+|X_{t^-}|}}dt+(1+|X_{t^-}|)^\frac{1}{4}dZ_t
\end{equation*}
where $(Z_t)$ is a purely discontinuous Lévy process having no
drift term with Lévy measure $\pi^{(\alpha)}$ (defined in the
preceding example). One checks that Proposition \ref{vitpart}
applies with $V(x)=1+|x|^2$, $a=3/4$, $r=1/4$ and every $p\in
(2,3)$. As in the preceding example, we test our procedure in the
cases $\alpha=1$ and $\alpha=5/3$. As the dynamical system is less
stable, the convergence is slower but we can observe the same
phenomenon.
\begin{figure}[h]
\begin{minipage}[c]{.5\linewidth}\centering
\includegraphics[width=5cm]{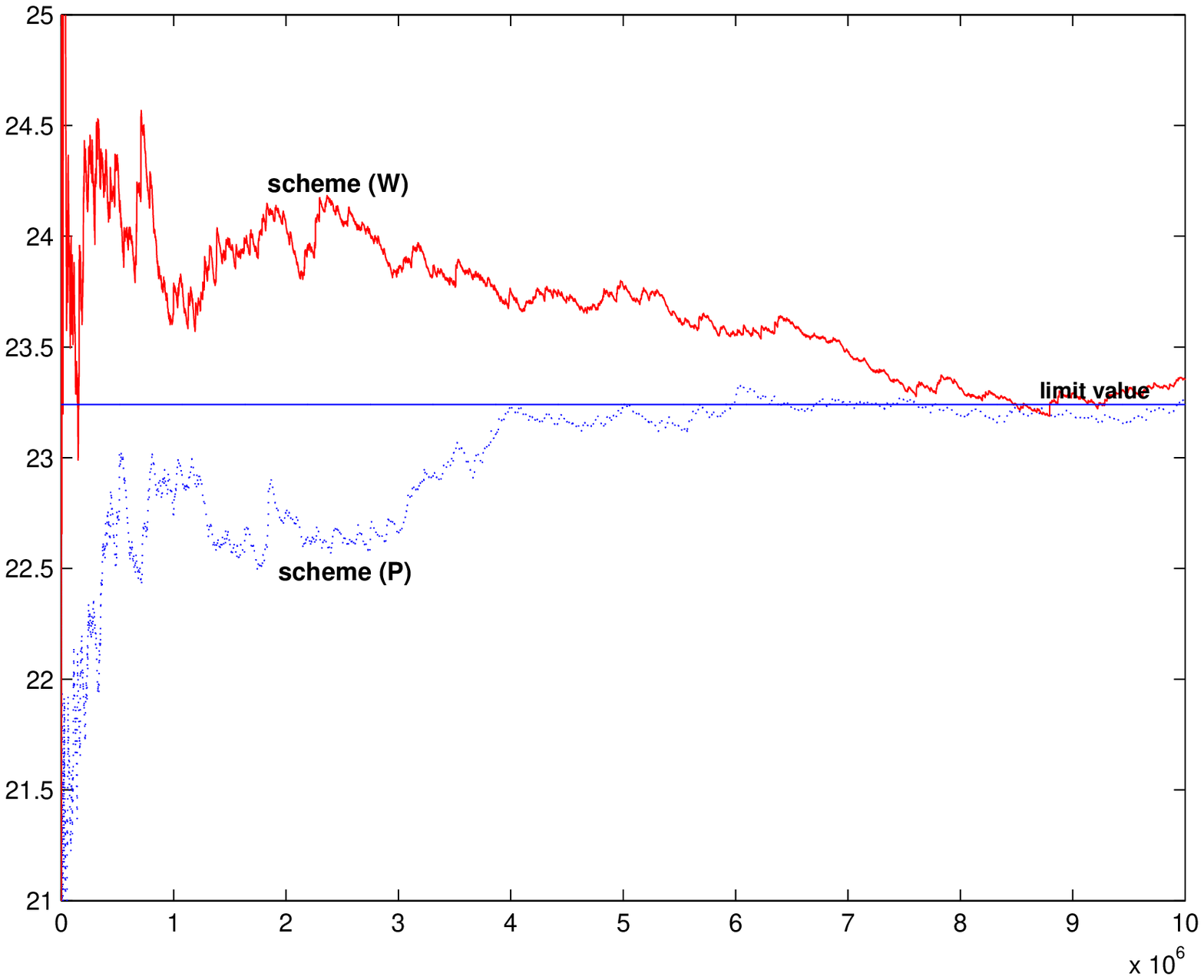}
\caption{$n\mapsto \bar{\nu}_n(\phi)$, $\alpha=1$}\label{figsim3}
\end{minipage}\hfill
\begin{minipage}[c]{.5\linewidth}\centering\label{figsim4}
\includegraphics[width=5cm]{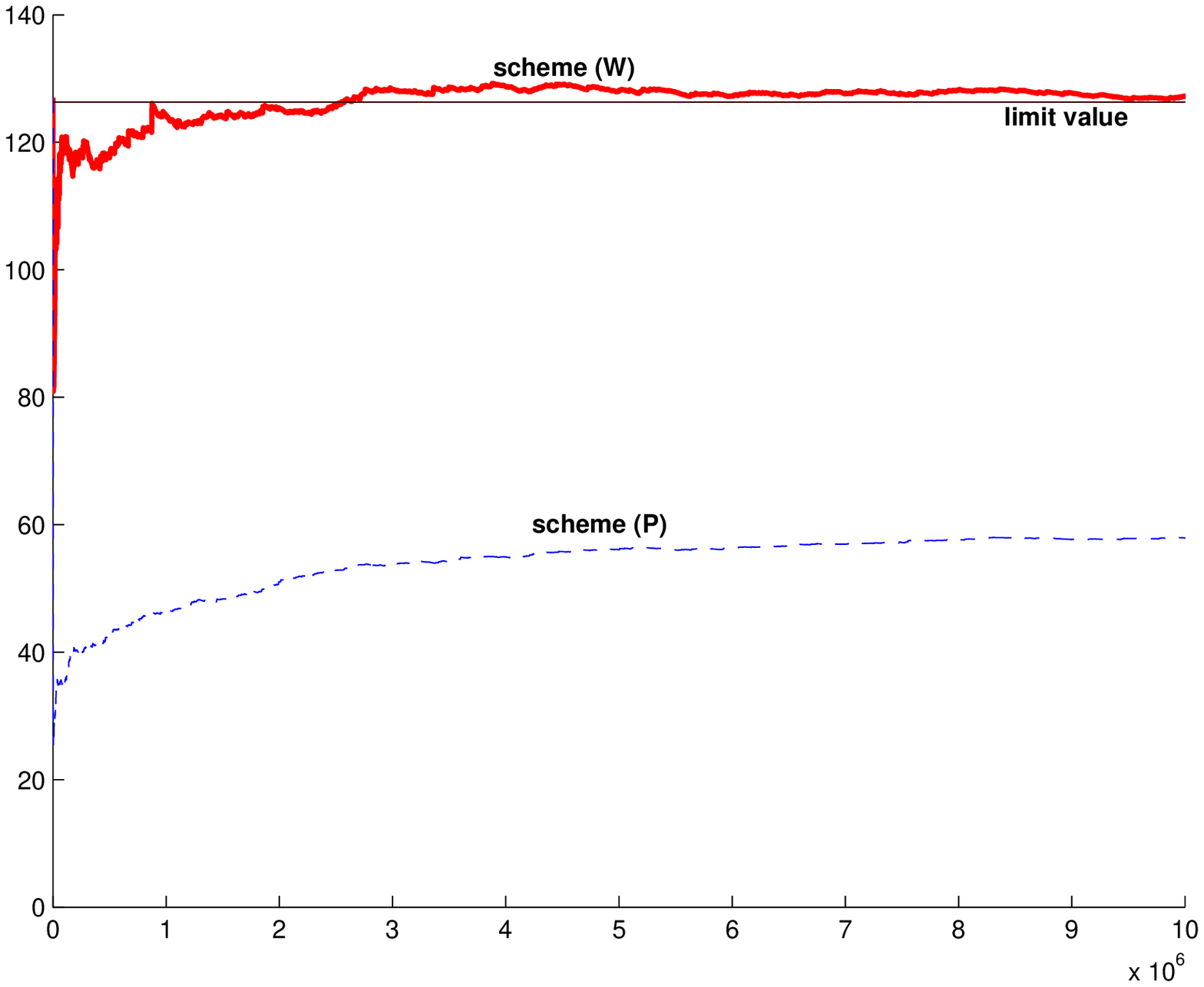}
\caption{$n\mapsto \bar{\nu}_n(\phi)$, $\alpha=5/3$}
\end{minipage}
\end{figure}


\begin{thebibliography}{99}
{\footnotesize
\bibitem[AsRo01]{bib10}
Asmussen S., Rosinski J. (2001), Approximations of small jumps of
Levy processes with a view towards simulation, J. of Applied
Probability, 38, pp. 482-493.
\bibitem[BaMiRe01]{bib18}
Barndorff-Nielsen O., Mikosh T., Resnick S. (2001), \textit{L\'evy
Processes and Applications}, Birkauser.
\bibitem[CoRo05]{cohen}
Cohen S., Rosinski J. (2005), Gaussian approximation of
multivariate Lévy processes with applications to simulation of
tempered and operator stable processes, Preprint.
\bibitem[EtKu86]{bib4}
Ethier S., Kurtz T. (1986), \textit{Markov processes,
characterization and convergence}, Wiley series in probability and
mathematical statistics, Wiley, New York.
\bibitem[HaHe80]{hall}
Hall P., Heyde C. (1980),  \textit{Martingale Limit Theory and its
Application}, Academic Press.
\bibitem[LaPa02]{bib2}
Lamberton D., Pag\`es G. (2002), Recursive computation of the
invariant distribution of a diffusion   \textit{Bernoulli}, 8, pp. 367-405.
\bibitem[LaPa03]{bib3}
Lamberton D., Pag\`es G. (2003), Recursive computation of the
invariant distribution of a diffusion: the case of a weakly mean
reverting drift, \textit{Stochastics and Dynamics}, 4,
pp. 435-451.
\bibitem[Lem05]{bib25}
Lemaire V. (2005), An adaptative scheme for the approximation of
dissipative systems, \textit{Stochastic Process. Appl.}, to
appear.
\bibitem[Lem06]{lemaire}
Lemaire V. (2005), \textit{Estimation numérique de la mesure
invariante d'une diffusion}, Phd Thesis,
http://tel.ccsd.cnrs.fr/tel-00011281.
\bibitem[Pag01]{bib14}
Pagès G. (2001), Sur quelques algorithmes récursifs pour les
probabilités numériques, \textit{ESAIM : Probab. Statis.}, 5,
pp. 141-170.
\bibitem[Pan05]{panloup}
Panloup F. (2005), Recursive computation of the invariant measure of a
stochastic differential equation driven by a Lévy process,
http://hal.ccsd.cnrs//ccsd-00009273, Preprint.
\bibitem[Pan06]{bib26} Panloup F., Phd Thesis, Université Paris VI, in progress.
\bibitem[PaVe01]{Parver1}
Pardoux E., Veretennikov A. (2001), On Poisson equation and diffusion approximation I, \textit{Ann. Probab.}, 29, pp. 1061-1085.
\bibitem[PaVe03]{Parver2}
Pardoux E., Veretennikov A. (2003), On Poisson equation and diffusion appoximation II, \textit{Ann. Probab}, 31 , pp. 1166-1192.
\bibitem[PaVe06]{Parver3}
Pardoux E., Veretennikov A. (2006),  On the Poisson equation and diffusion approximation III,
to appear in \textit{Ann. Probab}.}
\end{thebibliography}
\end{document}